\documentclass{amsart}
\usepackage{amsthm,graphicx}
\usepackage{amssymb}
\usepackage{colortbl}
\usepackage{tikz}

\usepackage[colorlinks=true]{hyperref}
\usepackage{textcomp}

\counterwithin{equation}{section}
\theoremstyle{definition}
 \newtheorem{theorem}{Theorem}
\newtheorem*{theorem*}{Theorem}
 \newtheorem{proposition}[equation]{Proposition}
 \newtheorem{definition}[equation]{Definition}
 \newtheorem{remark}[equation]{Remark}
 \newtheorem{lemma}[equation]{Lemma}
 \newtheorem{corollary}[equation]{Corollary}
 \newtheorem{example}[equation]{Example}

\newcommand\RR{\mathbb R}
\newcommand\QQ{\mathbb Q}
\newcommand\ZZ{\mathbb Z}
\newcommand\NN{\mathbb N}

\newcommand\A{\mathcal A}

\newcommand\state[1]{\langle #1\rangle}

\newcommand\e{\varepsilon}

\newcommand{\Add}{\mathrm{Add}}

\newcommand{\edge}{\mathrm{edge}}
\newcommand{\vertex}{\mathrm{vertex}}
\newcommand{\node}{\mathrm{node}}

\newcommand{\smooth}{\mathrm{Smooth}}
\newcommand{\p}{{\bf p}}
\newcommand{\q}{{\bf q}}
\newcommand{\h}{\text{\small h}}
\newcommand{\C}{\text{c}}

\numberwithin{equation}{section}

\begin{document}

\date{\today}
\title{Tropical curves in sandpile models}
\author[N. Kalinin, M. Shkolnikov]{Nikita Kalinin, Mikhail Shkolnikov}\thanks{Nikita Kalinin, Guangdong Technion Israel Institute of Technology (GTIIT),
Technion-Israel Institute of Technology,  ORCID: 0000-0002-1613-5175, nikaanspb@gmail.com. \\
Mikhail Shkolnikov, Institute of Mathematics and Informatics, Bulgarian Academy of Sciences, corresponding author,  m.shkolnikov@math.bas.bg, research is supported by the Simons Foundation International grant no. 992227, IMI-BAS}

\address{Guangdong Technion Israel Institute of Technology (GTIIT),
241 Daxue Road, Shantou, Guangdong Province 515603, P.R. China,Technion-Israel Institute of Technology, Haifa, 32000, Haifa District
Israel}

\address{Institute of Mathematics and Informatics, Bulgarian Academy of Sciences, Acad. G. Bonchev Str., Bl. 8, 1113 Sofia, Bulgaria.}

\begin{abstract}
A sandpile is a cellular automaton on a graph that evolves by the following toppling rule: if the number of grains at a vertex is at least its valency, then this vertex sends one grain to each of its neighbors.

In the study of pattern formation in sandpiles on large subgraphs of the standard square lattice, S. Caracciolo, G. Paoletti, and A. Sportiello experimentally observed that the result of the relaxation of a small perturbation of the maximal stable state contains a clear visible thin balanced graph formed by its deviation (less than maximum) set. Such graphs are known as tropical curves.

During the early stage of our research, we have noticed that these tropical curves are approximately scale-invariant, that is the deviation set mimics an extremal tropical curve depending on the domain on the plane and the positions of the perturbation points, but not on the mesh of the lattice.

In this paper, we rigorously formulate these two facts in the form of a scaling limit theorem and prove it. We rely on the theory of tropical analytic series, which is used to describe the global features of the sandpile dynamic, and on the theory of smoothings of discrete superharmonic functions, which handles local questions.
\end{abstract}
\maketitle



\begin{figure}[h!]
\includegraphics[scale=0.23]{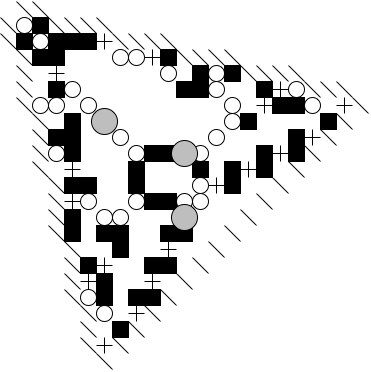}
\includegraphics[scale=0.23]{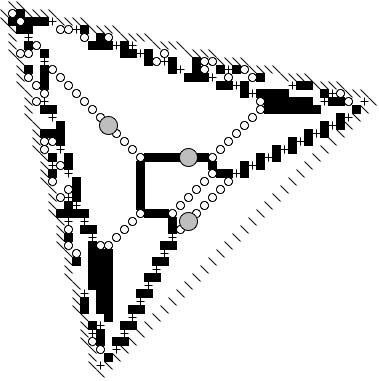}
\includegraphics[scale=0.23]{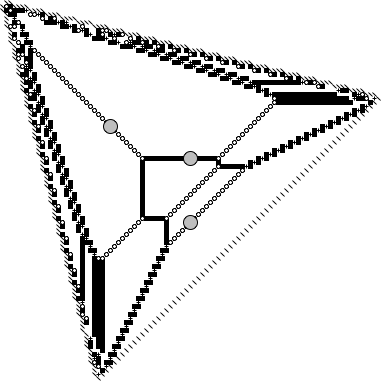}

\caption{A thin balanced graph appears as a deviation set of a sandpile. See Example~\ref{ex_evidence} for details. White corresponds to three grains, black to one, circles to two, crosses to zero, and skew lines are the boundary vertices (sinks). Grey rounds represent the positions of added grains.}
\label{fig_evidence}
\end{figure}


\section{Introduction}
{\bf History.} This article aims to establish a rigorous bridge between planar tropical curves \cite{mikh} and sandpiles on a square grid \cite{BTW, Dhar}. The first are often presented as combinatorial analogs of algebraic (and in our context, also analytic) curves, belonging to the field of tropical geometry where a polynomial (or series) is a piece-wise linear function with integral slopes. The second are most frequently affiliated with self-organized criticality, a concept in natural sciences, reflecting a critical behavior of the system similar to that near a phase transition, but without being finely tuned.

The discovery of such a relation between the two realms is not ours. It was predicted in \cite{firstsand}, and in the most common pictures of sandpiles, such as in the result of dropping many grains at the origin on the whole lattice or the identity element of the sandpile group for various domains \cite{melchionna2020sandpile,pegden2017stability}, discrete versions of tropical curves, i.e. thin balanced graphs made of deviations from the background, are clearly visible. However, the formal statement for their presence in these pictures is yet unknown. Our project started with the advice of G. Mikhalkin to have a look at the experimental works of S. Caracciolo, G. Paoletti, and A. Sportiello \cite{firstsand,book}, where such tropical curves appear in the vicinity of the maximal stable state.

To this observation, we add that the picture is approximately scale-invariant. Mathematically, this is expressed in the form of a scaling limit theorem, in its weaker form announced in \cite{announce}, and presented in the next section. The core of the proof is phenomenological in its essence. Locally, it relies on a careful analysis (established in \cite{us_solitons}) of the interplay between waves \cite{ivashkevich1994waves}, which are infinitesimal portions of sandpile avalanches, and solitons and triads, of which the discretized tropical curves are made. Global behavior is handled at the level of tropical series and ``grain'' operators $G_p$ acting on them \cite{us_series}, which are tropical counterparts of adding a grain at a point $p,$ relaxing, and removing a grain from $p$. 

Current work has inspired the definition of a continuous analog of the sandpile model, where a state is given by a finite collection of points $P$ in the interior of a convex domain together with a tropical series $F$ vanishing at the boundary of this domain. Such a state is called stable if $F$ is not smooth along $P$. For an unstable state, its relaxation consists of applying the grain operators $G_p$ for points $p$ from $P$, each contracting a face of the corner locus of $F$ until it passes through $p$. We don't know if such tropical relaxations always terminate after a finite number of steps (an instance of it is proven in the one-dimensional case \cite{shkolnikov2023relaxation}), but we know that it converges to some series $G_P F,$ independent of a particular relaxation, provided that all $G_p$ for $p\in P$ appear in it infinitely many times. The tropical sandpile model was studied numerically \cite{sandcomputation,stattrop}, in particular, it was shown that, similar to the usual sandpile model \cite{BTW}, it demonstrates a power-law for the distribution of avalanches, providing the first continuous example of self-organized criticality.

{\bf Results.} This paper is dedicated to proving that, in dimension two, the classical sandpile model in the near-maximal regime converges to the tropical sandpile model. More specifically, for a convex domain $\Omega$ and a collection of points $P\subset\Omega^\circ,$ the tropical analytic curve defined by the series $G_P 0_\Omega,$ where $0_\Omega$ on $\Omega$ is a tropical series equal to zero identically, approximates the locus of points where the maximal stable state on a lattice discretization of $\Omega$ loses grains after being perturbed at roundings of points in $P,$ provided that the lattice mesh is sufficiently small. Tropical analytic curves arising in such limits are interesting on their own, for example for $\Omega$ equal to a disk and $P$ consisting of a single point at its center, one gets the tropical caustic of this disk, see Fig. \ref{fig_disk}. Tropical caustic (the set of critical points of tropical wave front \cite{mikhalkin2023wave}), the concept which was revealed by the current work, is a complete invariant of convex domains, thus one can derive formulas for natural quantities such as area, perimeter, etc, in terms of their moduli (see \cite{pi_short,easter}). A somewhat surprising corollary of the naturality, i.e. its coordinate-free definition, of the tropical sandpile model is the emergence of modular invariance of the sandpile model in this regime for small lattice mesh. Tropical curves arising in this limit are extremal among all curves passing through the perturbation points -- they minimize action and symplectic area -- solving a tropical version of Steiner problem \cite{us_series}.

One crucial idea of the main theorem proof is that the relaxations in classical and tropical sandpile models are performed in parallel, i.e. each application of $G_p$ in the limit corresponds to sending some number of waves from a rounding of $p$ at the discrete level, provided that the sandpile state is made solely of solitons going along some at most nodal tropical curve, and waves do not reach the boundary. Another crucial idea is that by sending a slightly smaller number of waves each time, we ensure that we stay within this framework of solitons and at most nodal curves. Therefore, if we skillfully prepare the initial state by which we replace the maximal stable state, the partial relaxation goes smoothly and in a very controlled way. This partial relaxation gives a good approximation for the original perturbation of the maximal stable state and gives a lower bound on the toppling function of the latter. The upper bound, on the other hand, is easy to derive, it is achieved by simply rescaling and discretizing the tropical series $G_P 0_\Omega,$ and its validity follows from the least action principle. These two bounds are shown to be as close as it is requested when the mesh of the lattice tends to zero.

Our hope is that this geometric approach can be extended outside of the current regime. In particular, the interplay of waves and solitons that we exploit in the present article is just the ground level of the variety of other phenomena observed in the formation of sandpile patterns and their transformations. The next big step would be to explain how exactly the movement of tropical curves is responsible for the mass exchange and translation of quadratic patterns, as it is perceived, for example, in harmonic dynamics movies (see \cite{lang2019harmonic}), or in the scaling process of the identity element of the sandpile group, or in the growth of a single source sandpile. Overall, we expect that the result, and more importantly, the techniques presented in this paper are just the beginning of an exploration that will lead to a complete understanding of sandpiles on lattices.

\section{Main result}
\label{sec_main}
This section aims to formally state our scaling limit theorem in its most geometric form. Unfortunately, separated from its proof, the statement alone loses much of the acquired understanding of the sandpile model in the regime of a small perturbation of the maximal stable state on a large convex domain of the square lattice. Thus, a reader who is willing to gain some feeling of this regime and grasp why the theorem is true is invited to read the next section first. 
 

\subsection{Graphs in admissible domains}
\begin{definition}[\cite{us_series}]
\label{def_omegaadmissible} 
A convex closed subset $\Omega\subset\RR^2$ with non-empty interior $\Omega^\circ$ is said to be {\it admissible} if it does not contain a line with an irrational slope.
\end{definition}

Here is a complete list of convex closed but not admissible domains:
\begin{itemize}
\item  $\Omega$ is a subset of a line (i.e. $\Omega^\circ$ is empty);
\item $\Omega$ is $\RR^2$;
\item $\Omega$ is a half-plane with the boundary of an irrational slope;
\item $\Omega$ is a strip between two parallel lines of irrational slope.
\end{itemize} 

From now on we always suppose that $\Omega$ is an {\bf admissible} subset of $\RR^2$.

\begin{definition}
\label{def_lattice}
Consider the lattice $\h \ZZ^2=\{(i\h,j\h)|i,j\in\ZZ\}$ with the mesh $\h>0$, $\h\ZZ^2$ is seen as a set of vertices of a four-valent graph: we connect by edges the pairs $z,z'\in \h \ZZ^2$ of points with $|z-z'|=\h$ and denote the relation of being adjacent by $z\sim z'$. Let $\Omega_\h = \h\ZZ^2\cap\Omega^\circ$ and let $\partial\Omega_\h$ be the set of vertices of $\Omega_\h$ which have a neighbor vertex outside $\Omega^\circ$. 
\end{definition}

\subsection{Sandpiles}

A {\it state} $\phi$ of the sandpile model on $\Omega_\h$ is a function $$\phi:\Omega_\h\to \ZZ_{\geq 0}$$ on {vertices} of $\Omega_\h$. We interpret $\phi(z)$ as the number of grains of sand in $z\in\h\ZZ^2$. A vertex $z_0\in \Omega_\h\setminus \partial\Omega_\h$ is called {\it unstable} if $\phi(z_0)\geq 4$, and in this case $z_0$ can {\it topple}, producing a new state $\phi':\Omega_\h\to \ZZ_{\geq 0}$ by the following local rule changing only the value at $z_0$ and its neighbors:
\[
\phi'(z) = 
\begin{cases}
\phi(z)-4 & \text{if $z=z_0$;}\\
\phi(z)+1 & \text{if $z\sim z_0$;}\\
\phi(z)& \text{ otherwise.}\\
\end{cases}
\]
Note that we prohibit making topplings at the vertices in
$\partial\Omega_\h$ (or, equivalently, we think of $\partial\Omega_\h$ as sinks).
A ${\it relaxation}$ is doing topplings at unstable vertices while it is possible: if $\Omega_\h$ is finite, then any relaxation eventually terminates. For infinite $\Omega_\h$ for our purposes it is enough to consider so-called {\it locally-finite relaxations}, i.e. the relaxation where each vertex topples a finite number of times, see details in Appendix A of \cite{us_solitons}, this theory is quite parallel to the finite case.  

We denote by $\phi^\circ$ the result of a relaxation of $\phi$. It is a classical fact  \cite{Dhar} that $\phi^\circ$ does not depend on a particular relaxation and is uniquely determined by $\phi$. 
For a survey about sandpiles see \cite{MR3662912,MR2246566,redig,jarai2014sandpile} and references therein.  

Let $H_\phi:\Omega_\h\to \ZZ$ be the {\it toppling function} (also referred to as odometer in some other texts) of a state $\phi$, i.e. the function counting the number of topplings during a relaxation of $\phi$. Then,  \begin{equation}\label{eqn_relaxtoppling}\phi^\circ = \phi+\Delta (H_\phi), \end{equation} where the laplacian $\Delta F$ of a function is defined as

$$(\Delta F)(z) = -4F(z)+\sum_{z'\sim z}F(z').$$

On an infinite graph, for a state $\phi$ admitting a locally-finite relaxation, the toppling function $H_\phi$ is well-defined, and the formula \eqref{eqn_relaxtoppling} holds \cite{us_solitons}.

A function $F$ is harmonic (resp., superharmonic) on $A\subset \h\ZZ^2$ if $\Delta F$ = 0  (resp., $\Delta F\leq 0$) at every point in $A$. 


\subsection{Tropical series}

\begin{definition}
\label{def_tropicalanalitycal} 
An $\Omega$-{\it tropical series} is a continuous function 
$f:\Omega\to\RR_{\geq 0}$, $f|_{\partial\Omega}= 0$,  such that there exists $\A \subset \ZZ^2$ and $a_{ij}\in \RR$ for each $(i,j)\in \A$, and at each point  $(x,y)\in \Omega$ we have
\begin{equation}
\label{eq_series}
f(x,y)=\inf\limits_{(i,j)\in\A}(a_{ij}+ix+jy).
\end{equation} 
An $\Omega$-{\it
tropical curve} $C(f)$ on $\Omega^\circ$ is the corner locus of an $\Omega$-tropical series $f$, i.e. the set of points where $f$ is not smooth. 
\end{definition}

As we proved in \cite{us_series}, for each $\Omega$-tropical series $f$, for each $z\in\Omega^\circ$, $f$ is given by the minimum of a finite set of monomials $a_{ij}+ix+jy$ for $(x,y)$ in a small neighborhood of $z$. Therefore, locally, $C(f)$ is a graph with straight edges.

Note that such $f$, as a function from $\Omega$ to $\RR$, has many presentations in the form \eqref{eq_series}.  We suppose that $\A$ is chosen to be maximal by inclusion and the coefficients $a_{ij}$ are as minimal as possible. We call this presentation {\it the canonical form} of a tropical series. For each $\Omega$-tropical series, there exists a unique canonical form of it. For example, a canonical series form of a function $0_\Omega$ on $\Omega$ that is identically zero is $$0_\Omega(x,y)=\inf\limits_{(i,j)\in\A_\Omega}(c_{ij}^\Omega+ix+jy),$$ where $c_{ij}^\Omega$ is defined as $-\inf_{(x,y)\in\Omega}(ix+jy)$ for all $(i,j)\in\mathbb{Z}^2$ and $$\A_\Omega=\{(i,j):c_{ij}^\Omega\neq\infty\}.$$ Note that $(0,0)\in \A_\Omega$ and $c_{00}=0$. The property of being admissible for convex closed $\Omega$ with non-empty interior is equivalent to $\A_\Omega$ being strictly bigger than the one element set containing $(0,0).$ Therefore, for an admissible $\Omega$ we may define a non-trivial $\Omega$-tropical series $\rho_\Omega$ given (in a non-canonical form) simply by omitting the $(0,0)$-term in the canonical form of $0_\Omega.$ Morally, the function $\rho_\Omega$ measures a tropical distance to the boundary of $\Omega.$ The level sets of $\rho_\Omega$ are studied in \cite{mikhalkin2023wave}.  

\begin{definition}
\label{def_fOmegaP}
Let $\p_1,\dots \p_n\in \Omega^\circ$ be different points, $P=\{\p_1,\dots,\p_n\}$. We denote by $f_{\Omega,P}$ the pointwise minimum among all $\Omega$-tropical
series that are not smooth at all the points $\p_1,\dots,\p_n$.
\end{definition}

Observe that this definition requires that there is at least one $\Omega$-tropical series not smooth at all points $\p_1,\dots,\p_n$ -- such a series can be constructed with the help of the function $\rho_\Omega:$
$$\min(\rho_\Omega(z),\rho_\Omega(\p_1))+\dots+\min(\rho_\Omega(z),\rho_\Omega(\p_n)).$$ Here we implicitly use the fact that the set of $\Omega$-tropical series is naturally equipped with two operations -- the pointwise minimum -- playing the role of tropical addition, and the pointwise addition -- playing the role of tropical multiplication. In particular, the tropical curve defined by a usual pointwise sum of two series is a union of corresponding curves. 

\begin{lemma}[\cite{us_series}]
\label{lemma_minimalseries}
The function $f_{\Omega,P}$ is an $\Omega$-tropical series and $P\subset C(f_{\Omega,P}).$
\end{lemma}


\subsection{Main theorem}

\begin{definition}
\label{def_deviation}
Denote by $\state{3}$ the {\it maximal stable state}, the state which has exactly three grains at every vertex of $\Omega_\h$. 
 For a state $\phi$ on $\Omega_\h$, its {\it deviation locus} $D(\phi)$ is defined as  $$D(\phi)=\{z\in \Omega_\h | \phi(z)\ne 3 \}.$$
\end{definition}

S. Caracciolo, G. Paoletti, and A. Sportiello examined the result of a relaxation of a small perturbation (by adding grains to several points) of $\state{3}$ on $\Omega\cap \ZZ^2$. They observed \cite{firstsand,CPS} that the deviation locus of the relaxed state looks like a balanced graph (see Figure~\ref{fig_evidence}), also known as a tropical curve \cite{mikh}. This relation between sandpiles and tropical geometry was also mentioned in \cite{dhar2009pattern,dhar2013sandpile}. We add to this that these tropical curves a approximately scale invariant and rigorously formulate it as follows.


\begin{definition}
We say that $\p^\h\in \h\ZZ^2$ is {\bf a} rounding of a point $\p\in\RR^2$ with respect to
$\h\ZZ^2$ if the distance between $\p$ and $\p^\h$ is at most $\h$.
\end{definition}

Let $P$ be a finite subset of $\Omega^\circ$ and $P^\h = \{\p^\h|\p\in P\}$ be a set of proper roundings (Definition~\ref{def_proper}) of points in $P$ with
respect to the lattice $\h\ZZ^2$. Consider the state $\phi_\h$ of a sandpile on $\Omega_\h$ defined as
\begin{equation}
\label{eq_phih}
\phi_\h
= \state{3} + \sum_{\p\in P}\delta_{\p^\h}.
\end{equation}

Our main result is the following theorem announced in \cite{announce} in a weaker form. 
\begin{theorem}
\label{th_main}
The family of deviation sets $D(\phi_\h^\circ)$ converges (by Hausdorff, on
compact sets in $\Omega^\circ$) to the tropical curve $C(f_{\Omega,P})$ as $\h\to 0$. 
\end{theorem}

If $P$ is a one-element set, i.e. there is only one perturbation point, the theorem holds with any rounding. In the setup of computer experiments, when $\Omega$ is a lattice polygon, $P$ is a subset of $\mathbb{Z}^2$ and $\h^{-1}$ is an integer, the proper roundings can be chosen tautologically see Proposition \ref{prop_roundinglattice}, i.e. $ \p^\h=\p$ for all $\p\in P$ -- this is precisely the weaker form announced in \cite{announce}. This theorem has other versions that are also proven in the present paper that enhance the limit with a structure finer than that of merely a planar graph. For example, we establish the pointwise convergence of the corresponding toppling functions to $f_{\Omega,P}$ and show how to restore the multiplicities on the edges of $C(f_{\Omega,P})$ in yet another way through a certain weak-$^*$ limit. Multiplicities on edges (though, most of the time equal to $1$) are essential to formulating the geometric defining property of a planar tropical curve, i.e. {\it the balancing condition}: for every vertex, the sum of outgoing primitive vectors along the adjacent edges is zero when counted with the corresponding multiplicities.

\begin{example}
\label{ex_evidence}
Let $\Omega$ be a triangle given by three lines $$x-y=0, 4x+y=30,x+4y=120.$$ Let $\p_1,\p_2, \p_3$ be the points $(7,22),(12,20),(12,16),n=3$.  Figure~\ref{fig_evidence} shows the results of the relaxation  of $\phi_\h=\state{3}+\sum_{i=1}^n\delta_{\p_i}$ for $\h = 1/N$ where $N= 1,2,4$.  Pictures like that first appeared in \cite{firstsand}. It was observed by the authors of \cite{firstsand} that the deviation locus in Figure~\ref{fig_evidence} is balanced.
\end{example}

The ambiguity with roundings is justified as follows. The correspondence $$(\p_1,\p_2,\dots,\p_n)\to f_{\Omega,P}$$ is continuous for a generic set $P$ of points, and, in this case, Theorem~\ref{th_main} holds for any rounding $P^\h$ of $P$. But if $P$ belongs to the discriminant set\footnote{In our case this set is a subset of several hyperplanes (thus it has zero measure) in $\Omega\times \Omega\times \Omega\times \dots\times \Omega$ ($|P|$ times) and divides it into several chambers of full dimension.} of configurations, then $\phi_\h^\circ$ does not
depend in any sense continuously on the points where we drop additional sand;
the susceptibility of a sandpile is very dramatic. For different roundings of
$\p_1,\dots,\p_n$ with respect to $\h\ZZ^2$ we {\bf can} obtain drastically different pictures of $\phi_\h^\circ$ (e.g. see Figure 5.4 in \cite{misha_thesis}). Thus, for $P$ in the discriminant, we rather prove that there {\bf exist} so-called proper roundings $P^\h$ such that the above convergence result takes place. In fact, the proper roundings are the roundings for which the ``asymptotically minimal'' number of toppling occurs when relaxing $\phi_\h$.


The proof of Theorem~\ref{th_main} can be found in Section~\ref{sec_proof}. It is based on the following two technical tools, described in detail in subsequent sections. The first is local combinatorial analysis of discrete superharmonic function, developed earlier in \cite{us_solitons}: we prove that each rational direction $(p,q)\in\ZZ^2$ is represented by a pattern, so-called soliton or $(p,q)$-string, and a soliton moves changeless when applying sandpile waves (a relaxation can always be decomposed into waves). The second is global dynamics of these patterns can be very well controlled by means of tropical series theory, \cite{us_series}. We carve an approximation of any admissible domain $\Omega$ by rational slope polygons using level sets of the expected tropical series in the limit and construct tropical polynomials on these polygons which give us lower bounds for the relaxation of $\phi_\h$, and these lower bounds converge to the upper bound, achieved through the least action principle, as $\h\to 0$. 


\section{A non-technical overview}
\subsection{The simplest example: single perturbation point on a standard square}
It is instructive to go to the basics and consider in detail the example of a single perturbation point on a square domain. To be completely precise, take $\Omega=[0,10]^2$ and $\p=(3,4).$ To represent states on $\Omega_1=\Omega^\circ\cap\mathbb{Z}^2$ we draw it as a $9\times 9$ chekered square and write numbers in the square sites. We don't write $3$ since it is expected to be the most frequent number of grains at a vertex, as on Fig. \ref{fig_beginrelax34} presenting the beginning of a particular relaxation for $\state{3}+\delta_\p$, the perturbation of the maximal stable state at $\p.$ Every further step of this relaxation is defined by doing topplings (removing four grains from a vertex and giving one to each neighbor) at all unstable vertices (the numbers of grains at which are greater than three)  of the previous state. In order to derive the result of the relaxation, we apply the crucial trick of decomposing relaxations into waves.

\begin{figure}
\begin{tikzpicture}[scale=.3, every node/.style={scale=.7}]
\draw (1, 1) grid (10, 10);

\node at (3.5,4.5) {4};

\begin{scope}[xshift=300]
\draw (1, 1) grid (10, 10);

\node at (3.5,4.5) {0};

\node at (2.5,4.5) {4};
\node at (3.5,5.5) {4};
\node at (3.5,3.5) {4};
\node at (4.5,4.5) {4};
\end{scope}

\begin{scope}[xshift=600]
\draw (1, 1) grid (10, 10);

\node at (3.5,4.5) {4};

\node at (2.5,4.5) {0};
\node at (3.5,5.5) {0};
\node at (3.5,3.5) {0};
\node at (4.5,4.5) {0};

\node at (1.5,4.5) {4};
\node at (5.5,4.5) {4};
\node at (3.5,6.5) {4};
\node at (3.5,2.5) {4};

\node at (2.5,3.5) {5};
\node at (2.5,5.5) {5};
\node at (4.5,3.5) {5};
\node at (4.5,5.5) {5};

\end{scope}

\begin{scope}[xshift=900]
\draw (1, 1) grid (10, 10);

\node at (3.5,4.5) {0};

\node at (2.5,4.5) {4};
\node at (3.5,5.5) {4};
\node at (3.5,3.5) {4};
\node at (4.5,4.5) {4};

\node at (1.5,4.5) {0};
\node at (5.5,4.5) {0};
\node at (3.5,6.5) {0};
\node at (3.5,2.5) {0};

\node at (2.5,3.5) {1};
\node at (2.5,5.5) {1};
\node at (4.5,3.5) {1};
\node at (4.5,5.5) {1};

\node at (6.5,4.5) {4};
\node at (3.5,7.5) {4};
\node at (3.5,1.5) {4};

\node at (1.5,5.5) {5};
\node at (2.5,6.5) {5};
\node at (4.5,6.5) {5};
\node at (5.5,5.5) {5};
\node at (5.5,3.5) {5};
\node at (4.5,2.5) {5};
\node at (2.5,2.5) {5};
\node at (1.5,3.5) {5};

\end{scope}
\end{tikzpicture}
\caption{The first few steps of a relaxation on a $9\times 9$ square where we add a single grain to the all-$3$-state, empty sites represnt the value $3$.}
\label{fig_beginrelax34}
\end{figure}

Sending a wave from a vertex $\p$ is a partially defined operator $W_\p$ on the space of stable states on the graph which is a composition of the (forced) toppling at $\p$ and a subsequent relaxation, see Fig. \ref{fig_9by9onewave}.  The mentioned toppling is called forced since it is applied to a stable vertex which results in a negative value at this vertex. Therefore, one needs to extend the notion of relaxation to such states with possibly negative values (it is totally straightforward). Note that for $W_\p\phi$ to be an honest non-negative-valued state the vertex $\p$ must have a neighbor $\q$ such that $\phi(\q)$ is the maximal stable value, i.e. $3$ in case of square lattice subgraphs. Waves are relatively easy to understand -- essentially, applying a wave $W_\p$ is doing one toppling at each vertex of the cluster of threes containing $\p$ -- this property follows from an easy fact that $W_\p\phi=W_\q\phi$ if $\p$ and $\q$ are adjacent vertices such that $\phi(\p)=\phi(\q)=3.$   Be aware that doing just this may result in a non-stable state, i.e. there are more topplings to perform in some other vertices.

\begin{figure}
\begin{tikzpicture}[scale=.3, every node/.style={scale=.7}]

\draw (1, 1) grid (10, 10);

\node at (3.5,4.5) {-1};

\node at (2.5,4.5) {4};
\node at (3.5,5.5) {4};
\node at (3.5,3.5) {4};
\node at (4.5,4.5) {4};

\begin{scope}[xshift=300]
\draw (1, 1) grid (10, 10);

\node at (2.5,4.5) {0};
\node at (3.5,5.5) {0};
\node at (3.5,3.5) {0};
\node at (4.5,4.5) {0};

\node at (1.5,4.5) {4};
\node at (5.5,4.5) {4};
\node at (3.5,6.5) {4};
\node at (3.5,2.5) {4};

\node at (2.5,3.5) {5};
\node at (2.5,5.5) {5};
\node at (4.5,3.5) {5};
\node at (4.5,5.5) {5};

\end{scope}

\begin{scope}[xshift=600]
\draw (1, 1) grid (10, 10);

\node at (1.5,4.5) {0};
\node at (5.5,4.5) {0};
\node at (3.5,6.5) {0};
\node at (3.5,2.5) {0};

\node at (2.5,3.5) {1};
\node at (2.5,5.5) {1};
\node at (4.5,3.5) {1};
\node at (4.5,5.5) {1};

\node at (6.5,4.5) {4};
\node at (3.5,7.5) {4};
\node at (3.5,1.5) {4};

\node at (1.5,5.5) {5};
\node at (2.5,6.5) {5};
\node at (4.5,6.5) {5};
\node at (5.5,5.5) {5};
\node at (5.5,3.5) {5};
\node at (4.5,2.5) {5};
\node at (2.5,2.5) {5};
\node at (1.5,3.5) {5};

\end{scope}

\begin{scope}[xshift=900]
\draw (1, 1) grid (10, 10);

\node at (1.5,4.5) {2};

\node at (6.5,4.5) {0};
\node at (3.5,7.5) {0};
\node at (3.5,1.5) {0};

\node at (1.5,5.5) {1};
\node at (2.5,6.5) {1};
\node at (4.5,6.5) {1};
\node at (5.5,5.5) {1};
\node at (5.5,3.5) {1};
\node at (4.5,2.5) {1};
\node at (2.5,2.5) {1};
\node at (1.5,3.5) {1};

\node at (3.5,8.5) {4};
\node at (7.5,4.5) {4};

\node at (1.5,6.5) {5};
\node at (2.5,7.5) {5};

\node at (4.5,7.5) {5};
\node at (5.5,6.5) {5};
\node at (6.5,5.5) {5};

\node at (4.5,1.5) {5};
\node at (5.5,2.5) {5};
\node at (6.5,3.5) {5};

\node at (1.5,2.5) {5};
\node at (2.5,1.5) {5};
\end{scope}

\begin{scope}[yshift=-300]
\draw (1, 1) grid (10, 10);

\node at (1.5,4.5) {2};
\node at (1.5,5.5) {2};
\node at (1.5,3.5) {2};
\node at (3.5,1.5) {2};

\node at (3.5,8.5) {0};
\node at (7.5,4.5) {0};

\node at (1.5,6.5) {1};
\node at (2.5,7.5) {1};

\node at (4.5,7.5) {1};
\node at (5.5,6.5) {1};
\node at (6.5,5.5) {1};

\node at (4.5,1.5) {1};
\node at (5.5,2.5) {1};
\node at (6.5,3.5) {1};

\node at (1.5,2.5) {1};
\node at (2.5,1.5) {1};

\node at (3.5,9.5) {4};
\node at (8.5,4.5) {4};

\node at (1.5,1.5) {5};
\node at (1.5,7.5) {5};
\node at (2.5,8.5) {5};
\node at (4.5,8.5) {5};
\node at (5.5,7.5) {5};
\node at (6.5,6.5) {5};
\node at (7.5,5.5) {5};
\node at (7.5,3.5) {5};
\node at (6.5,2.5) {5};
\node at (5.5,1.5) {5};
\end{scope}

\begin{scope}[yshift=-300,xshift=300]
\draw (1, 1) grid (10, 10);

\node at (1.5,4.5) {2};
\node at (1.5,5.5) {2};
\node at (1.5,3.5) {2};
\node at (3.5,1.5) {2};
\node at (1.5,6.5) {2};
\node at (1.5,2.5) {2};

\node at (4.5,1.5) {2};
\node at (2.5,1.5) {2};

\node at (3.5,9.5) {0};
\node at (8.5,4.5) {0};

\node at (1.5,1.5) {1};
\node at (1.5,7.5) {1};
\node at (2.5,8.5) {1};
\node at (4.5,8.5) {1};
\node at (5.5,7.5) {1};
\node at (6.5,6.5) {1};
\node at (7.5,5.5) {1};
\node at (7.5,3.5) {1};
\node at (6.5,2.5) {1};
\node at (5.5,1.5) {1};

\node at (9.5,4.5) {4};

\node at (1.5,8.5) {5};
\node at (2.5,9.5) {5};
\node at (4.5,9.5) {5};
\node at (5.5,8.5) {5};
\node at (6.5,7.5) {5};
\node at (7.5,6.5) {5};
\node at (8.5,5.5) {5};
\node at (8.5,3.5) {5};
\node at (7.5,2.5) {5};
\node at (6.5,1.5) {5};
\end{scope}

\begin{scope}[yshift=-300,xshift=600]
\node at (5.5,6) {\Large five more steps};
\node at (5.3,5) {\Large $\longrightarrow$};
\end{scope}

\begin{scope}[yshift=-300,xshift=900]
\draw (1, 1) grid (10, 10);

\node at (1.5,1.5) {1};
\node at (1.5,9.5) {1};
\node at (9.5,1.5) {1};
\node at (9.5,9.5) {1};

\node at (1.5,2.5) {2};
\node at (1.5,3.5) {2};
\node at (1.5,4.5) {2};
\node at (1.5,5.5) {2};
\node at (1.5,6.5) {2};
\node at (1.5,7.5) {2};
\node at (1.5,8.5) {2};
\node at (9.5,2.5) {2};
\node at (9.5,3.5) {2};
\node at (9.5,4.5) {2};
\node at (9.5,5.5) {2};
\node at (9.5,6.5) {2};
\node at (9.5,7.5) {2};
\node at (9.5,8.5) {2};
\node at (2.5,1.5) {2};
\node at (3.5,1.5) {2};
\node at (4.5,1.5) {2};
\node at (5.5,1.5) {2};
\node at (6.5,1.5) {2};
\node at (7.5,1.5) {2};
\node at (8.5,1.5) {2};
\node at (2.5,9.5) {2};
\node at (3.5,9.5) {2};
\node at (4.5,9.5) {2};
\node at (5.5,9.5) {2};
\node at (6.5,9.5) {2};
\node at (7.5,9.5) {2};
\node at (8.5,9.5) {2};

\end{scope}

\end{tikzpicture}
\caption{Propagation of a wave on all-$3$-background initiated by forcing a toppling at a site, procuring the value $-1$ at this site at the initial step. As before, empty sites represent the value $3$. When the wave hits the boundary some grains are lost (washed away) and a less-than-three value ($2$ along the edges and $1$ at the corners) is left behind at the next step of the propagation. Note that the result is independent of the source of the wave and can be defined simply by doing exactly one toppling at every site.}
\label{fig_9by9onewave}
\end{figure}

Any relaxation can be decomposed into waves. More specifically, for a state of the form $\state{3}+\delta_\p$ there exist $m\in\mathbb{N}$ such that $$(\state{3}+\delta_\p)^\circ=\delta_\p+W_\p^m\state{3}.$$ This $m$ is the minimal number such that $W_\p^m\state{3}$ is less than $3$ at $\p$. For example, if $\p$ is adjacent to the boundary of the domain it is enough to apply a single wave and then drop a grain at $\p.$ Otherwise, one needs to apply more waves, see Fig. \ref{fig_9by9waves2345}. For example, for $\p=(3,4)$ three waves are necessary, and for $\p=(5,5)$ -- five waves. In general, the number of waves needed to complete a relaxation of $\state{3}+\delta_p$ can be seen as the lattice distance from $\p$ to the boundary.  

\begin{figure}
\begin{tikzpicture}[scale=.3, every node/.style={scale=.7}]
\draw (1, 1) grid (10, 10);

\node at (1.5,1.5) {1};
\node at (2.5,2.5) {1};

\node at (1.5,9.5) {1};
\node at (2.5,8.5) {1};

\node at (9.5,1.5) {1};
\node at (8.5,2.5) {1};

\node at (9.5,9.5) {1};
\node at (8.5,8.5) {1};

\node at (2.5,3.5) {2};
\node at (2.5,4.5) {2};
\node at (2.5,5.5) {2};
\node at (2.5,6.5) {2};
\node at (2.5,7.5) {2};

\node at (8.5,3.5) {2};
\node at (8.5,4.5) {2};
\node at (8.5,5.5) {2};
\node at (8.5,6.5) {2};
\node at (8.5,7.5) {2};

\node at (3.5,2.5) {2};
\node at (4.5,2.5) {2};
\node at (5.5,2.5) {2};
\node at (6.5,2.5) {2};
\node at (7.5,2.5) {2};

\node at (3.5,8.5) {2};
\node at (4.5,8.5) {2};
\node at (5.5,8.5) {2};
\node at (6.5,8.5) {2};
\node at (7.5,8.5) {2};

\begin{scope}[xshift=300]
\draw (1, 1) grid (10, 10);

\node at (1.5,1.5) {1};
\node at (2.5,2.5) {1};
\node at (3.5,3.5) {1};

\node at (1.5,9.5) {1};
\node at (2.5,8.5) {1};
\node at (3.5,7.5) {1};

\node at (9.5,1.5) {1};
\node at (8.5,2.5) {1};
\node at (7.5,3.5) {1};

\node at (9.5,9.5) {1};
\node at (8.5,8.5) {1};
\node at (7.5,7.5) {1};

\node at (3.5,4.5) {2};
\node at (3.5,5.5) {2};
\node at (3.5,6.5) {2};

\node at (7.5,4.5) {2};
\node at (7.5,5.5) {2};
\node at (7.5,6.5) {2};

\node at (4.5,3.5) {2};
\node at (5.5,3.5) {2};
\node at (6.5,3.5) {2};

\node at (4.5,7.5) {2};
\node at (5.5,7.5) {2};
\node at (6.5,7.5) {2};
\end{scope}

\begin{scope}[xshift=600]
\draw (1, 1) grid (10, 10);

\node at (1.5,1.5) {1};
\node at (2.5,2.5) {1};
\node at (3.5,3.5) {1};
\node at (4.5,4.5) {1};

\node at (1.5,9.5) {1};
\node at (2.5,8.5) {1};
\node at (3.5,7.5) {1};
\node at (4.5,6.5) {1};

\node at (9.5,1.5) {1};
\node at (8.5,2.5) {1};
\node at (7.5,3.5) {1};
\node at (6.5,4.5) {1};

\node at (9.5,9.5) {1};
\node at (8.5,8.5) {1};
\node at (7.5,7.5) {1};
\node at (6.5,6.5) {1};

\node at (4.5,5.5) {2};

\node at (6.5,5.5) {2};

\node at (5.5,4.5) {2};

\node at (5.5,6.5) {2};
\end{scope}

\begin{scope}[xshift=900]
\draw (1, 1) grid (10, 10);

\node at (1.5,1.5) {1};
\node at (2.5,2.5) {1};
\node at (3.5,3.5) {1};
\node at (4.5,4.5) {1};

\node at (1.5,9.5) {1};
\node at (2.5,8.5) {1};
\node at (3.5,7.5) {1};
\node at (4.5,6.5) {1};

\node at (9.5,1.5) {1};
\node at (8.5,2.5) {1};
\node at (7.5,3.5) {1};
\node at (6.5,4.5) {1};

\node at (9.5,9.5) {1};
\node at (8.5,8.5) {1};
\node at (7.5,7.5) {1};
\node at (6.5,6.5) {1};

\node at (5.5,5.5) {-1};
\end{scope}

\end{tikzpicture}
\caption{The results of sending two, three, four, and five waves from the center of the square.}
\label{fig_9by9waves2345}
\end{figure}

The approximate scale invariance of the resulting pictures is apparent, see Fig \ref{fig_relaxone34}. The totality of sites where the value deviates from $3$ is called the deviation locus. Our main result is that deviation loci have scaling limits which are tropical curves passing through perturbation points. In fact, these tropical curves are clearly visible at a discrete level provided that the mesh $\h$ of the lattice is sufficiently small with respect to the number of perturbation points. These discrete analogs of tropical curves are made of soliton patterns. In this basic example, we have encountered solitons for four slopes, i.e. $0$ -- horizontal, $\infty$ -- vertical (both are made of the value $2$), and $\pm 1$ -- two diagonals (made of the value $1$). Patterns for other slopes are more intricate, we will see another one (up to symmetries of the square tiling) while considering perturbations of the maximal stable state on a square at two points and  on Fig. \ref{fig_triads}.

\begin{figure}
\begin{tikzpicture}[scale=.38, every node/.style={scale=.9}]
\draw[very thin, gray!20](0,0) grid (10,10);

\draw(0,0)--(10,0)--(10,10)--(0,10)--(0,0);

\node at (1,1) {1};
\node at (2,2) {1};
\node at (3,3) {1};

\node at (1,9) {1};
\node at (2,8) {1};
\node at (3,7) {1};

\node at (9,1) {1};
\node at (8,2) {1};
\node at (7,3) {1};

\node at (9,9) {1};
\node at (8,8) {1};
\node at (7,7) {1};

\node at (3,5) {2};
\node at (3,6) {2};

\node at (7,4) {2};
\node at (7,5) {2};
\node at (7,6) {2};

\node at (4,3) {2};
\node at (5,3) {2};
\node at (6,3) {2};

\node at (4,7) {2};
\node at (5,7) {2};
\node at (6,7) {2};

\begin{scope}[xshift=330, every node/.style={scale=.7}]
\draw[step=0.5, very thin, gray!20](0,0) grid (10,10);

\draw(0,0)--(10,0)--(10,10)--(0,10)--(0,0);

\node at (0.5,0.5) {1};
\node at (1,1) {1};
\node at (1.5,1.5) {1};
\node at (2,2) {1};
\node at (2.5,2.5) {1};
\node at (3,3) {1};

\node at (0.5,9.5) {1};
\node at (1,9) {1};
\node at (1.5,8.5) {1};
\node at (2,8) {1};
\node at (2.5,7.5) {1};
\node at (3,7) {1};

\node at (9.5,0.5) {1};
\node at (9,1) {1};
\node at (8.5,1.5) {1};
\node at (8,2) {1};
\node at (7.5,2.5) {1};
\node at (7,3) {1};

\node at (9.5,9.5) {1};
\node at (9,9) {1};
\node at (8.5,8.5) {1};
\node at (8,8) {1};
\node at (7.5,7.5) {1};
\node at (7,7) {1};

\node at (3,3.5) {2};

\node at (3,4.5) {2};
\node at (3,5) {2};
\node at (3,5.5) {2};
\node at (3,6) {2};
\node at (3,6.5) {2};

\node at (7,3.5) {2};
\node at (7,4) {2};
\node at (7,4.5) {2};
\node at (7,5) {2};
\node at (7,5.5) {2};
\node at (7,6) {2};
\node at (7,6.5) {2};

\node at (3.5,3) {2};
\node at (4,3) {2};
\node at (4.5,3) {2};
\node at (5,3) {2};
\node at (5.5,3) {2};
\node at (6,3) {2};
\node at (6.5,3) {2};

\node at (3.5,7) {2};
\node at (4,7) {2};
\node at (4.5,7) {2};
\node at (5,7) {2};
\node at (5.5,7) {2};
\node at (6,7) {2};
\node at (6.5,7) {2};

\end{scope}

\begin{scope}[xshift=660]
\draw(0,0)--(10,0)--(10,10)--(0,10)--(0,0);
\draw(0,0)--(3,3)--(3,7)--(0,10);
\draw(10,10)--(7,7)--(3,7);
\draw(10,0)--(7,3)--(3,3);
\draw(7,3)--(7,7);

\node at (3,4) {$\bullet$};

\end{scope}

\end{tikzpicture}

\caption{Left: the result of the relaxation started on Fig. \ref{fig_beginrelax34}, compare it with the second picture on Fig. \ref{fig_9by9waves2345}, i.e. the result of applying three waves. Note the difference in the presentation of squares and numbers: on the previous pictures, for the sake of visibility, we have drawn the numbers inside squares that we called sites; now, however, the numbers are written at vertices of light gray squares, the exact relation between these and former squares is the duality between two square tilings, i.e. faces of one are vertices of another.  Center: rescaling by factor $2$ of the domain on the left, i.e. the set of vertices is $(0,10)^2\cap{1\over 2}\mathbb{Z}^2$ (instead of $(0,10)^2\cap\mathbb{Z}^2$) and the grain is dropped at the same vertex $(3,4)$. Right: the scaling limit of deviation loci given by a tropical curve passing through the perturbation point $(3,4)$.}
\label{fig_relaxone34}
\end{figure}

An important thing to note about this example is that only the first wave interacts with the boundary of the square, creating four edges made of solitons, and that these edges, and especially their conjunctions, move in a predictable way under the action of subsequent waves. Similar phenomena take place in general, but deriving it requires considerable effort. To be more precise, for a general rational slope, one needs multiple waves to create a soliton pattern going along a boundary line segment with this slope. In the case of corners of the boundary formed by two segments, the picture is more delicate, being dependent on the modularity of this corner. If the corner is unimodular, i.e. the primitive vectors going in the directions of its sides form a basis of the lattice, the behavior under the wave action (after applying first a sufficient number of waves) is the same as for the standard quadrant corner, i.e. there are two soliton edges parallel to the sides of the corner that move in the direction of the source of waves and leave behind a single deviation edge made of a soliton pattern going along the trajectory of the vertex at the intersection of these two edges (where this direction is the sum of primitive vectors in the directions of the sides of the corner).

To complete the description of the most basic example we would like to write explicitly the corresponding toppling function. Let $\Omega=[0,1]^2$ and $\p\in\Omega_\h=\Omega^\circ\cap \h\mathbb{Z}^2$ for some $\h>0.$ Denote by $H$ the toppling function of the state $\state{3}+\delta_\p.$ As we observed before, the number of waves needed to relax this state is the $\h\mathbb{Z}^2$-lattice distance $\rho_\h(\p)$ from $\p$ to $\partial\Omega.$ Moreover, the number of topplings made at a vertex $\q$ during the relaxation is $\rho_\h(\q)$ if $\rho_\h(\q)<\rho_\h(\p),$ i.e. only this number of waves sent from $\p$ reach $\q$, and, otherwise equal to $\rho_\h(\p)$ since all waves reach $q.$ This gives $$H(\q)=\min(\rho_\h(\q),\rho_\h(\p)),$$ where $$\rho_\h(x,y)=\h^{-1}\min(x,y,1-x,1-y).$$ We see that, at least in this example, $\h H$ is a restriction of a tropical polynomial to the $\h$-lattice points of $\Omega.$

\subsection{Directions of generalization}
After understanding completely the case of a single perturbation point on a standard square, one may ask themselves: ``How does it generalize?'' Here are some directions:
\begin{enumerate}  
\item Replacing the square with more general domains. The technically simplest case would be that of unimodular convex lattice polygons, i.e. those for which at every vertex a pair of primitive vectors in the directions of adjacent sides is a basis of the lattice,  the next step would be to pass to convex polygons with rational slope sides, and further step is of general compact convex domains; dropping compactness is not hard, but we need to exclude those domains that contain a line of irrational slope for the purposes of applying tropical series theory (also, at the sandpile level different scaling is needed in order to derive a non-trivial limit) and we also need to extend the notion of finite relaxation by a locally finite one. Dropping convexity, on the other hand, is beyond our current technique at the tropical level, also sandpiles behave differently on non-convex domains (see Fig. \ref{fig_nonconvex}).
\begin{figure}
\label{fig_nonconvex}
\includegraphics[width=0.84\textwidth]{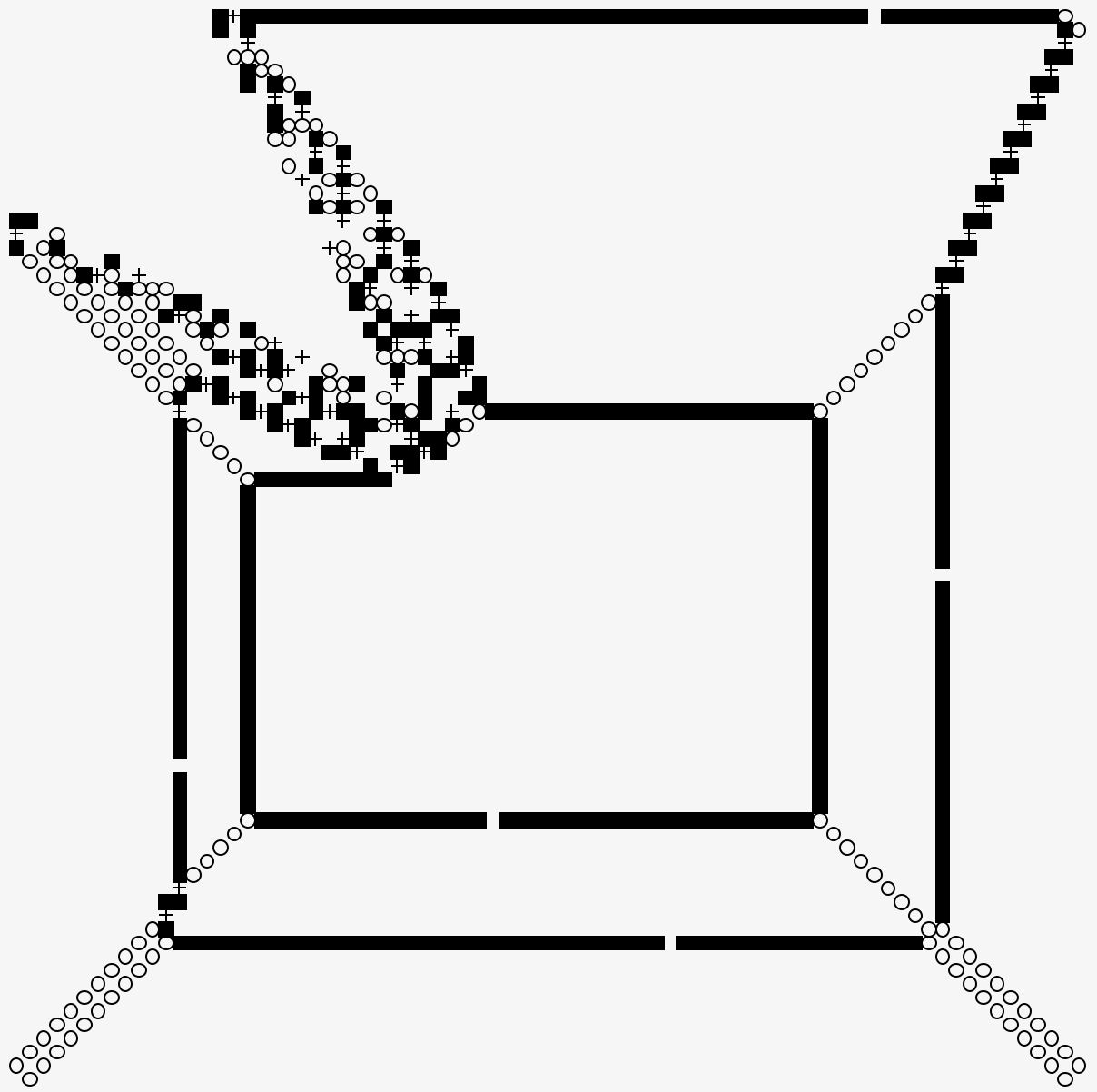}
\caption{An example of how a perturbation of the maximal stable state on a non-convex domain may look like. }
\end{figure}

\item Increasing number of perturbation points, in this article we deal with an arbitrary finite subset of the interior of the domain, for which at a rounding of each of its points we add a single grain to the maximal stable state. In contrast with a single perturbation point, the choice of such roundings is not always arbitrary, i.e. for non-generic configurations of points, there are preferred roundings for which the limit of toppling function is the predicted tropical series (otherwise the limit over subsequences would still exist, but it may be greater). Another difficulty we need to face comes from the current gap of understanding of stabilization for relaxations in the tropical sandpile model, however, it is not that grave to obstruct completely the proof of the scaling limit theorems, but makes it slightly more involved. Both rounding and stabilization technicalities disappear in the most common case of computer simulations, i.e. when the domain is a lattice polygon, all perturbation points belong to the lattice and the inverse of the scaling factor is integral. What we don't allow in the current setup is to drop more than one grain in a perturbation point, this is due to two reasons. At the sandpile level, we would need to use multi-solitons, that are not yet developed. At the tropical level, we don't understand completely what it means for a tropical curve to pass through a given point with a given multiplicity (and this multiplicity unfortunately will not be related to the number of grains dropped to this point).
\item The maximal stable state could be replaced with a more general heavy background, i.e. it can be taking values less than three at some points. (if it is too light, for example consisting of the value two everywhere, there is nothing to study since no avalanches will be produced by adding a few extra grains). A physical condition for such a background, in order for our techniques to be still applicable, would be that waves of topplings can propagate through it from any region to another. Of course, a precise formal result in this generality would be even messier in terms of roundings, since a rounding of a perturbation point could be dropped on a less than three value, in which case no avalanche is produced.
\item Our theory can be extended to other regular graphs on the Euclidean plane, such as provided by triangular or hexagonal tilings. For higher order regular tilings (of hyperbolic plane), however, the situation is very different -- we tested it only for the heptagonal tilling \cite{kalinin2015sandpiles}, but the general conceptual difference with the zero curvature case is the following: perturbing the maximal stable state results in losing too much sand from the system, so we don't observe a linear deviation locus, since the number of vertices at the boundary of a big domain is proportional to the total number of vertices. On the square lattice, on the other hand, a big domain has only a linear number of boundary vertices with respect to the inverse of the scaling parameter as opposed to the quadratic total of the vertices in the domain -- this is why we expect to see a linear deviation locus for perturbing the maximal stable state in our setup.
\item The first author has already started to develop the theory for extending our theorems to higher dimensions. This consists of the local analysis at the sandpile level \cite{Kalinin_2023} and the global analysis at the level of tropical geometry \cite{kalinin2021shrinking}.
\item One could also think of the plane as the simplest example of a tropical surface, it would be interesting to investigate sandpiles on other tropical varieties, where sending a wave would correspond to an infinitesimal deformation of a hypersurface.
\end{enumerate}

\subsection{Single perturbation point on a unimodular polygon}
Perhaps the most natural direction of generalizing the basic example is by considering domains different from the standard square. Below, we explain the case of a unimodular polygon as an underlying domain of the model. 

Let's start with the simplest admisible domain $\Omega$, i.e. a half-plane whose boundary line has a rational slope, and a point $\p\in\Omega^\circ.$ Assume for simplicity that the origin belongs to the boundary of $\Omega.$  One of the key results of \cite{us_solitons} is that after applying some number of waves from a rounding $\p_\h\in\Omega_\h$ (where $\h$ is sufficiently small) to the state $\state{3}$ on $\Omega_\h$ results in detaching a soliton pattern from the boundary. Applying more such waves moves the soliton away from the boundary with a minimal possible non-zero speed and doesn't change the values near the boundary. A relaxation of the state $\state{3}+\delta_{\p_\h}$ on $\Omega_\h$ is seen as applying waves until the value at $\p_\h$ is less than $3,$ i.e. the soliton has reached the point, and then adding a grain at $\p_\h.$ Therefore the toppling function of $\state{3}+\delta_{\p_\h}$ is approximately equal to $z\in\Omega_\h\mapsto \h^{-1}\min(\lambda(z),\lambda(\p_\h)),$ where $\lambda$ is the minimal support function of $\Omega$ with integer gradient.

Next, we assume that $\Omega$ is a unimodular lattice cone with an apex at the origin. Again we add a single grain at an $\h$-rounding $\p^\h\in\Omega_h$ of some point $\p\in\Omega^\circ$ and we want to compute (or rather estimate) the corresponding toppling function. Let $\lambda_1$ and $\lambda_2$ be the minimal support functions of the cone $\Omega$ with integral gradients such that each vanishes on one of the boundary rays. Let $\rho\colon\Omega\rightarrow\mathbb{R}_{\geq 0}$ be defined as  $\rho(z)=\min(\lambda_1(z),\lambda_2(z)).$ We claim that the toppling function of $\state{3}+\delta_{\p_\h}$ on $\Omega_\h$ is approximated for small $\h>0$ by $\h^{-1}f_{\Omega,\p}|_{\Omega_\h},$ where $f_{\Omega,\p}\colon\Omega\rightarrow\mathbb{R}_{\geq 0}$ is given by $$f_{\Omega,\p}(z)=\min(\rho(z),\rho(\p)).$$

More specifically, since the function $(\h^{-1}f_{\Omega,p})|_{\Omega_\h}$ has a negative Laplacian at $\p_\h,$ it gives an upper bound on the toppling function by the least action principle \cite{FLP}. Technically (see section \ref{sec_proproundestimabove}), we choose the rounding $\p_\h\Omega_\h$ of $\p$ in order to have this property -- however, for one perturbation point, this is not necessary and $\p_\h$ can be chosen as an arbitrary $\h$-lattice point near $\p$ by the price of a slight modification of the estimate, i.e. by $\pm \h$ shift in the constant term $\rho(\p)$ of $f_{\Omega,\p}(z)=\min(\rho(z),\rho(\p))$.

The lower bound $(\h^{-1}(f_{\Omega,\p})|_{\Omega_\h}-\varepsilon)$ on the toppling function for any $\varepsilon>0$ and sufficiently small $\h>0,$ is achieved by using a restriction of a {\bf triad} to the cone. Triads are sandpile versions of smooth tropical curves having a single vertex and three (infinite) edges. Such a tropical curve is a corner locus of a tropical polynomial $F\colon\mathbb{R}^2\rightarrow\mathbb{R}$ of a form $$F(x,y)=\min(c_0,c_1+a_1x+b_1y,c_2+a_2x+b_2y),$$ where $a_1,b_1,a_2,b_2\in\mathbb{Z}$ and $a_1b_2-a_2b_1=\pm 1.$ Assuming in addition that $c_0,c_1,c_2$ are also integer, we may consider a sequence of functions $$F_0,F_1,F_2,\dots \colon \mathbb{Z}^2\rightarrow \mathbb{Z}$$ starting with $F_0=F|_\mathbb{Z}^2$ and satisfying the property that $F_n$ (called the $n$-smoothing) is the pointwise minimum of all superharmonic functions that are bounded from above by $F_0,$ by $F_0-n$ from below and differ from $F_0$ in some finite neighborhood of the corner locus $C(F).$ By the the main theorem of \cite{us_solitons} this sequence stabilizes at some function $\theta_F.$ A triad is defined by $F$ is $\state{3}+\Delta \theta_F.$

Construct an auxiliary stable state $\psi_{\h,\varepsilon}$ on $\Omega_h$ which is a restriction of the triad defined by $\h^{-1}\min(\rho,\varepsilon)$, i.e. it is made of two solitonic rays going along the boundary rays of the cone at the lattice distance $\varepsilon$ from them and a solitonic segment going from the apex of the cone to the junction of the solitonic rays. Triads behave in a predictable way under the action of waves (i.e. each wave increases by one the corresponding coefficient of the tropical polynomial and the triad itself moves towards the source of the wave with the minimal possible speed), therefore we know the toppling function of the state $\psi_{\h,\varepsilon}+\delta_{\p_\h}.$ This toppling function, which we estimate as $\h^{-1}(f_{\Omega,\p}|_{\Omega_\h}-\varepsilon),$ can be seen as a toppling function of a partial relaxation of $\state{3}+\delta_{\p_\h},$ where we withhold the excess grains along the deviation locus of $\psi_{\h,\varepsilon}.$  

The above argument can be repeated verbatim for a general unimodular polygon $\Omega$ with the only difference that now $\rho$ is defined as a minimum of primitive support functions for each side of $\Omega.$ This results in having a generally more complicated structure of the auxiliary state $\psi_{\h,\varepsilon}$ which now is potentially a combination of more than one triad. Technically speaking, we perform a sufficient smoothing of an extension of $\h^{-1}\min(\rho,\varepsilon)$ from $\Omega_\h$ to the lattice $\h\mathbb{Z}^2$ and then restrict it back to $\Omega_\h$.  Applying $n$ waves from $\p_\h$ to this state would result in the same order smoothing of $\h^{-1}\min(\rho,\varepsilon+n \h).$ In this process, some solitonic edges might contract resulting in singularities of the tropical curve (see Fig. \ref{fig_oneptpent}). However, these singularities, unless we are near a terminal state of the relaxation, are just simple nodes by a Proposition \ref{prop_atmostnodal}, for which the canonical smoothings are also well defined. Therefore, we retain perfect control in such relaxations and know that the corresponding toppling function is not much less than $\h^{-1}f_{\Omega,\p}.$ 

\begin{figure}
\begin{tikzpicture}[scale=.45, every node/.style={scale=.7}]

\draw[step=0.5, very thin, gray!20](0,0) grid (12,12);
\draw(0,0)--(0,12)--(10,12)--(12,10)--(12,0)--(0,0);

\node at (3,4) {$4$};

\node at (.5,.5) {1};
\node at (.5,1) {2};
\node at (.5,1.5) {2};
\node at (.5,2) {2};
\node at (.5,2.5) {2};
\node at (.5,3) {2};
\node at (.5,3.5) {2};
\node at (.5,4) {2};
\node at (.5,4.5) {2};
\node at (.5,5) {2};
\node at (.5,5.5) {2};
\node at (.5,6) {2};
\node at (.5,6.5) {2};
\node at (.5,7) {2};
\node at (.5,7.5) {2};
\node at (.5,8) {2};
\node at (.5,8.5) {2};
\node at (.5,9) {2};
\node at (.5,9.5) {2};
\node at (.5,10) {2};
\node at (.5,10.5) {2};
\node at (.5,11) {2};
\node at (.5,11.5) {1};

\node at (1,.5) {2};
\node at (1.5,.5) {2};
\node at (2,.5) {2};
\node at (2.5,.5) {2};
\node at (3,.5) {2};
\node at (3.5,.5) {2};
\node at (4,.5) {2};
\node at (4.5,.5) {2};
\node at (5,.5) {2};
\node at (5.5,.5) {2};
\node at (6,.5) {2};
\node at (6.5,.5) {2};
\node at (7,.5) {2};
\node at (7.5,.5) {2};
\node at (8,.5) {2};
\node at (8.5,.5) {2};
\node at (9,.5) {2};
\node at (9.5,.5) {2};
\node at (10,.5) {2};
\node at (10.5,.5) {2};
\node at (11,.5) {2};
\node at (11.5,.5) {1};

\node at (11.5,1) {2};
\node at (11.5,1.5) {2};
\node at (11.5,2) {2};
\node at (11.5,2.5) {2};
\node at (11.5,3) {2};
\node at (11.5,3.5) {2};
\node at (11.5,4) {2};
\node at (11.5,4.5) {2};
\node at (11.5,5) {2};
\node at (11.5,5.5) {2};
\node at (11.5,6) {2};
\node at (11.5,6.5) {2};
\node at (11.5,7) {2};
\node at (11.5,7.5) {2};
\node at (11.5,8) {2};
\node at (11.5,8.5) {2};
\node at (11.5,9) {2};
\node at (11.5,9.5) {2};
\node at (11.5,10) {1};

\node at (1,11.5) {2};
\node at (1.5,11.5) {2};
\node at (2,11.5) {2};
\node at (2.5,11.5) {2};
\node at (3,11.5) {2};
\node at (3.5,11.5) {2};
\node at (4,11.5) {2};
\node at (4.5,11.5) {2};
\node at (5,11.5) {2};
\node at (5.5,11.5) {2};
\node at (6,11.5) {2};
\node at (6.5,11.5) {2};
\node at (7,11.5) {2};
\node at (7.5,11.5) {2};
\node at (8,11.5) {2};
\node at (8.5,11.5) {2};
\node at (9,11.5) {2};
\node at (9.5,11.5) {2};
\node at (10,11.5) {1};

\node at (10.5,11) {1};
\node at (11,10.5) {1};

\begin{scope}[xshift=370]
\draw[step=0.5, very thin, gray!20](0,0) grid (12,12);
\draw(0,0)--(0,12)--(10,12)--(12,10)--(12,0)--(0,0);

\node at (3,4) {$4$};

\node at (.5,.5) {1};

\node at (1,1) {1};
\node at (1,1.5) {2};
\node at (1,2) {2};
\node at (1,2.5) {2};
\node at (1,3) {2};
\node at (1,3.5) {2};
\node at (1,4) {2};
\node at (1,4.5) {2};
\node at (1,5) {2};
\node at (1,5.5) {2};
\node at (1,6) {2};
\node at (1,6.5) {2};
\node at (1,7) {2};
\node at (1,7.5) {2};
\node at (1,8) {2};
\node at (1,8.5) {2};
\node at (1,9) {2};
\node at (1,9.5) {2};
\node at (1,10) {2};
\node at (1,10.5) {2};
\node at (1,11) {1};

\node at (.5,11.5) {1};

\node at (1.5,1) {2};
\node at (2,1) {2};
\node at (2.5,1) {2};
\node at (3,1) {2};
\node at (3.5,1) {2};
\node at (4,1) {2};
\node at (4.5,1) {2};
\node at (5,1) {2};
\node at (5.5,1) {2};
\node at (6,1) {2};
\node at (6.5,1) {2};
\node at (7,1) {2};
\node at (7.5,1) {2};
\node at (8,1) {2};
\node at (8.5,1) {2};
\node at (9,1) {2};
\node at (9.5,1) {2};
\node at (10,1) {2};
\node at (10.5,1) {2};
\node at (11,1) {1};

\node at (11.5,.5) {1};

\node at (11,1.5) {2};
\node at (11,2) {2};
\node at (11,2.5) {2};
\node at (11,3) {2};
\node at (11,3.5) {2};
\node at (11,4) {2};
\node at (11,4.5) {2};
\node at (11,5) {2};
\node at (11,5.5) {2};
\node at (11,6) {2};
\node at (11,6.5) {2};
\node at (11,7) {2};
\node at (11,7.5) {2};
\node at (11,8) {2};
\node at (11,8.5) {2};
\node at (11,9) {2};
\node at (11,9.5) {2};
\node at (11,10) {1};

\node at (11.5,10) {2};

\node at (1.5,11) {2};
\node at (2,11) {2};
\node at (2.5,11) {2};
\node at (3,11) {2};
\node at (3.5,11) {2};
\node at (4,11) {2};
\node at (4.5,11) {2};
\node at (5,11) {2};
\node at (5.5,11) {2};
\node at (6,11) {2};
\node at (6.5,11) {2};
\node at (7,11) {2};
\node at (7.5,11) {2};
\node at (8,11) {2};
\node at (8.5,11) {2};
\node at (9,11) {2};
\node at (9.5,11) {2};
\node at (10,11) {1};

\node at (10,11.5) {2};

\node at (10.5,10.5) {1};

\end{scope}

\begin{scope}[yshift=-370]
\draw[step=0.5, very thin, gray!20](0,0) grid (12,12);
\draw(0,0)--(0,12)--(10,12)--(12,10)--(12,0)--(0,0);

\node at (3,4) {$4$};

\node at (.5,.5) {1};

\node at (1,1) {1};

\node at (1.5,1.5) {1};
\node at (1.5,2) {2};
\node at (1.5,2.5) {2};
\node at (1.5,3) {2};
\node at (1.5,3.5) {2};
\node at (1.5,4) {2};
\node at (1.5,4.5) {2};
\node at (1.5,5) {2};
\node at (1.5,5.5) {2};
\node at (1.5,6) {2};
\node at (1.5,6.5) {2};
\node at (1.5,7) {2};
\node at (1.5,7.5) {2};
\node at (1.5,8) {2};
\node at (1.5,8.5) {2};
\node at (1.5,9) {2};
\node at (1.5,9.5) {2};
\node at (1.5,10) {2};
\node at (1.5,10.5) {1};

\node at (1,11) {1};

\node at (.5,11.5) {1};

\node at (2,1.5) {2};
\node at (2.5,1.5) {2};
\node at (3,1.5) {2};
\node at (3.5,1.5) {2};
\node at (4,1.5) {2};
\node at (4.5,1.5) {2};
\node at (5,1.5) {2};
\node at (5.5,1.5) {2};
\node at (6,1.5) {2};
\node at (6.5,1.5) {2};
\node at (7,1.5) {2};
\node at (7.5,1.5) {2};
\node at (8,1.5) {2};
\node at (8.5,1.5) {2};
\node at (9,1.5) {2};
\node at (9.5,1.5) {2};
\node at (10,1.5) {2};
\node at (10.5,1.5) {1};

\node at (11,1) {1};

\node at (11.5,.5) {1};

\node at (10.5,2) {2};
\node at (10.5,2.5) {2};
\node at (10.5,3) {2};
\node at (10.5,3.5) {2};
\node at (10.5,4) {2};
\node at (10.5,4.5) {2};
\node at (10.5,5) {2};
\node at (10.5,5.5) {2};
\node at (10.5,6) {2};
\node at (10.5,6.5) {2};
\node at (10.5,7) {2};
\node at (10.5,7.5) {2};
\node at (10.5,8) {2};
\node at (10.5,8.5) {2};
\node at (10.5,9) {2};
\node at (10.5,9.5) {2};

\node at (11,10) {2};

\node at (11.5,10) {2};

\node at (2,10.5) {2};
\node at (2.5,10.5) {2};
\node at (3,10.5) {2};
\node at (3.5,10.5) {2};
\node at (4,10.5) {2};
\node at (4.5,10.5) {2};
\node at (5,10.5) {2};
\node at (5.5,10.5) {2};
\node at (6,10.5) {2};
\node at (6.5,10.5) {2};
\node at (7,10.5) {2};
\node at (7.5,10.5) {2};
\node at (8,10.5) {2};
\node at (8.5,10.5) {2};
\node at (9,10.5) {2};
\node at (9.5,10.5) {2};

\node at (10,11) {2};

\node at (10,11.5) {2};

\node at (10,10.5) {1};
\node at (10.5,10) {1};

\end{scope}

\begin{scope}[yshift=-370, xshift=370]
\draw[step=0.5, very thin, gray!20](0,0) grid (12,12);
\draw(0,0)--(0,12)--(10,12)--(12,10)--(12,0)--(0,0);

\node at (3,4) {$4$};

\node at (.5,.5) {1};

\node at (1,1) {1};

\node at (1.5,1.5) {1};

\node at (2,2) {1};
\node at (2,2.5) {2};
\node at (2,3) {2};
\node at (2,3.5) {2};
\node at (2,4) {2};
\node at (2,4.5) {2};
\node at (2,5) {2};
\node at (2,5.5) {2};
\node at (2,6) {2};
\node at (2,6.5) {2};
\node at (2,7) {2};
\node at (2,7.5) {2};
\node at (2,8) {2};
\node at (2,8.5) {2};
\node at (2,9) {2};
\node at (2,9.5) {2};
\node at (2,10) {1};

\node at (1.5,10.5) {1};

\node at (1,11) {1};

\node at (.5,11.5) {1};

\node at (2.5,2) {2};
\node at (3,2) {2};
\node at (3.5,2) {2};
\node at (4,2) {2};
\node at (4.5,2) {2};
\node at (5,2) {2};
\node at (5.5,2) {2};
\node at (6,2) {2};
\node at (6.5,2) {2};
\node at (7,2) {2};
\node at (7.5,2) {2};
\node at (8,2) {2};
\node at (8.5,2) {2};
\node at (9,2) {2};
\node at (9.5,2) {2};
\node at (10,2) {1};

\node at (10.5,1.5) {1};

\node at (11,1) {1};

\node at (11.5,.5) {1};

\node at (10,2.5) {2};
\node at (10,3) {2};
\node at (10,3.5) {2};
\node at (10,4) {2};
\node at (10,4.5) {2};
\node at (10,5) {2};
\node at (10,5.5) {2};
\node at (10,6) {2};
\node at (10,6.5) {2};
\node at (10,7) {2};
\node at (10,7.5) {2};
\node at (10,8) {2};
\node at (10,8.5) {2};
\node at (10,9) {2};
\node at (10,9.5) {2};

\node at (11,10) {2};

\node at (11.5,10) {2};

\node at (2.5,10) {2};
\node at (3,10) {2};
\node at (3.5,10) {2};
\node at (4,10) {2};
\node at (4.5,10) {2};
\node at (5,10) {2};
\node at (5.5,10) {2};
\node at (6,10) {2};
\node at (6.5,10) {2};
\node at (7,10) {2};
\node at (7.5,10) {2};
\node at (8,10) {2};
\node at (8.5,10) {2};
\node at (9,10) {2};
\node at (9.5,10) {2};

\node at (10,11) {2};

\node at (10,11.5) {2};

\node at (10,10.5) {2};
\node at (10.5,10) {2};

\node at (10,10) {1};

\end{scope}

\begin{scope}[yshift=-740]
\draw[step=0.5, very thin, gray!20](0,0) grid (12,12);
\draw(0,0)--(0,12)--(10,12)--(12,10)--(12,0)--(0,0);

\node at (3,4) {$4$};

\node at (.5,.5) {1};

\node at (1,1) {1};

\node at (1.5,1.5) {1};

\node at (2,2) {1};

\node at (2.5,3) {2};
\node at (2.5,3.5) {2};
\node at (2.5,4) {2};
\node at (2.5,4.5) {2};
\node at (2.5,5) {2};
\node at (2.5,5.5) {2};
\node at (2.5,6) {2};
\node at (2.5,6.5) {2};
\node at (2.5,7) {2};
\node at (2.5,7.5) {2};
\node at (2.5,8) {2};
\node at (2.5,8.5) {2};
\node at (2.5,9) {2};

\node at (2,10) {1};

\node at (1.5,10.5) {1};

\node at (1,11) {1};

\node at (.5,11.5) {1};

\node at (2.5,2.5) {1};

\node at (3,2.5) {2};
\node at (3.5,2.5) {2};
\node at (4,2.5) {2};
\node at (4.5,2.5) {2};
\node at (5,2.5) {2};
\node at (5.5,2.5) {2};
\node at (6,2.5) {2};
\node at (6.5,2.5) {2};
\node at (7,2.5) {2};
\node at (7.5,2.5) {2};
\node at (8,2.5) {2};
\node at (8.5,2.5) {2};
\node at (9,2.5) {2};

\node at (10,2) {1};

\node at (10.5,1.5) {1};

\node at (11,1) {1};

\node at (11.5,.5) {1};

\node at (9.5,2.5) {1};
\node at (9.5,3) {2};
\node at (9.5,3.5) {2};
\node at (9.5,4) {2};
\node at (9.5,4.5) {2};
\node at (9.5,5) {2};
\node at (9.5,5.5) {2};
\node at (9.5,6) {2};
\node at (9.5,6.5) {2};
\node at (9.5,7) {2};
\node at (9.5,7.5) {2};
\node at (9.5,8) {2};
\node at (9.5,8.5) {2};
\node at (9.5,9) {2};
\node at (9.5,9.5) {1};

\node at (11,10) {2};

\node at (11.5,10) {2};

\node at (2.5,9.5) {1};
\node at (3,9.5) {2};
\node at (3.5,9.5) {2};
\node at (4,9.5) {2};
\node at (4.5,9.5) {2};
\node at (5,9.5) {2};
\node at (5.5,9.5) {2};
\node at (6,9.5) {2};
\node at (6.5,9.5) {2};
\node at (7,9.5) {2};
\node at (7.5,9.5) {2};
\node at (8,9.5) {2};
\node at (8.5,9.5) {2};
\node at (9,9.5) {2};

\node at (10,11) {2};

\node at (10,11.5) {2};

\node at (10,10.5) {2};
\node at (10.5,10) {2};

\node at (10,10) {1};

\end{scope}

\begin{scope}[yshift=-740, xshift=370]
\draw[step=0.5, very thin, gray!20](0,0) grid (12,12);
\draw(0,0)--(0,12)--(10,12)--(12,10)--(12,0)--(0,0);

\node at (.5,.5) {1};

\node at (1,1) {1};

\node at (1.5,1.5) {1};

\node at (2,2) {1};

\node at (3,3) {1};
\node at (3,3.5) {2};

\node at (3,4.5) {2};
\node at (3,5) {2};
\node at (3,5.5) {2};
\node at (3,6) {2};
\node at (3,6.5) {2};
\node at (3,7) {2};
\node at (3,7.5) {2};
\node at (3,8) {2};
\node at (3,8.5) {2};
\node at (3,9) {1};

\node at (2,10) {1};

\node at (1.5,10.5) {1};

\node at (1,11) {1};

\node at (.5,11.5) {1};

\node at (2.5,2.5) {1};

\node at (3.5,3) {2};
\node at (4,3) {2};
\node at (4.5,3) {2};
\node at (5,3) {2};
\node at (5.5,3) {2};
\node at (6,3) {2};
\node at (6.5,3) {2};
\node at (7,3) {2};
\node at (7.5,3) {2};
\node at (8,3) {2};
\node at (8.5,3) {2};
\node at (9,3) {1};

\node at (10,2) {1};

\node at (10.5,1.5) {1};

\node at (11,1) {1};

\node at (11.5,.5) {1};

\node at (9.5,2.5) {1};

\node at (9,3) {1};
\node at (9,3.5) {2};
\node at (9,4) {2};
\node at (9,4.5) {2};
\node at (9,5) {2};
\node at (9,5.5) {2};
\node at (9,6) {2};
\node at (9,6.5) {2};
\node at (9,7) {2};
\node at (9,7.5) {2};
\node at (9,8) {2};
\node at (9,8.5) {2};

\node at (9,9) {1};

\node at (9.5,9.5) {1};

\node at (11,10) {2};

\node at (11.5,10) {2};

\node at (2.5,9.5) {1};

\node at (3.5,9) {2};
\node at (4,9) {2};
\node at (4.5,9) {2};
\node at (5,9) {2};
\node at (5.5,9) {2};
\node at (6,9) {2};
\node at (6.5,9) {2};
\node at (7,9) {2};
\node at (7.5,9) {2};
\node at (8,9) {2};
\node at (8.5,9) {2};

\node at (10,11) {2};

\node at (10,11.5) {2};

\node at (10,10.5) {2};
\node at (10.5,10) {2};

\node at (10,10) {1};

\end{scope}

\end{tikzpicture}
\caption{Relaxation by sending waves for a single perturbation point on a unimodular pentagon. Observe that there is an intermediate step (after sending three waves), where the underlying tropical curve has a simple node. This is the only singularity type that we must work with.}
\label{fig_oneptpent}
\end{figure}

\subsection{Unimodularity dropped: single perturbation on rotated squares}
\begin{figure}\label{fig_rotatedsquares}
\includegraphics[width=0.84\textwidth]{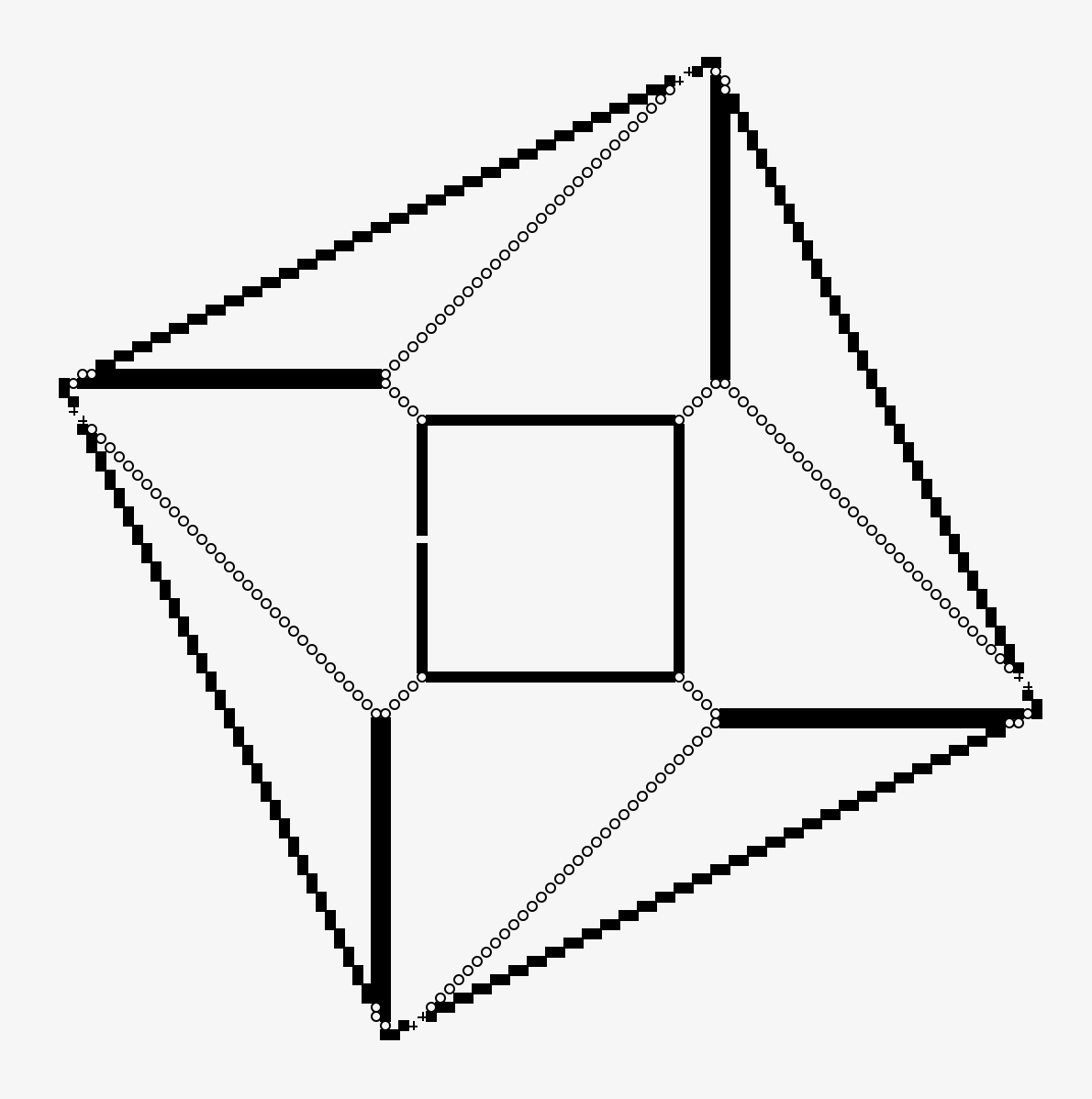}
\caption{An example for dropping a single grain on a rotated square. Observe that in contrast to the unimodular case, two edges are adjacent to each corner, one of which is double. In general, the number of edges going to a corner and their multiplicities is determined by a continued fraction expansion of a tropical cotangent of this corner. In particular, corners with exactly one adjacent edge correspond to $A_n$-corners (in terms of the singularity type of a symplectic toric surface having the polygon as its moment domain), and the multiplicity of this edge is $n+1.$ }
\end{figure}

If the unimodularity condition for a corner of a domain is dropped, the limiting tropical curve will have a more complicated incidence with the tip of the corner. Namely, it can have more edges ending at the tip of the corner, and/or some of these edges may have multiplicity greater than one. Observe for example Fig. 3.9 of \cite{misha_thesis} where such limiting curves are shown for various squares. At the discrete level, such a curve looks like in Fig. \ref{fig_rotatedsquares}. In terms of deriving results for such domains, it is not hard to reduce the case of a general polygon with rational slope sides to the case of unimodular polygons by cutting some corners. In toric geometry terminology, this corresponds to a well-known technique of resolution of singularities via blow-ups. 

\subsection{Tropical caustic of a disk}
\begin{figure}\label{fig_disk}
\includegraphics[width=0.48\textwidth]{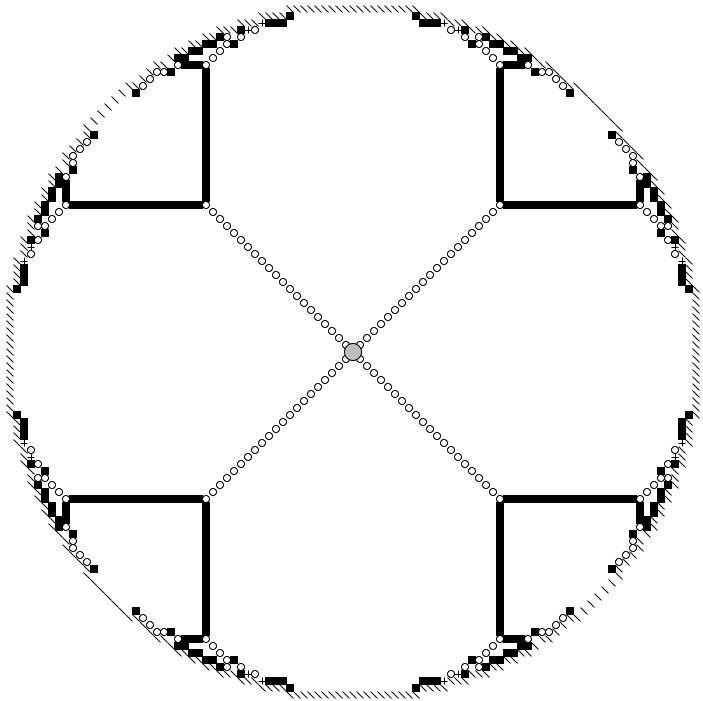}
\includegraphics[width=0.48\textwidth]{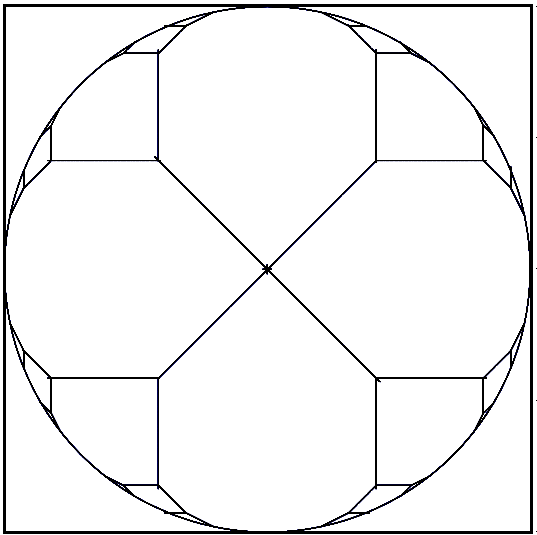}
\caption{Left: the result of adding a single grain at the center of a disk to the maximal stable state and relaxing. Right: tropical caustic of the disk.}
\end{figure}

The first example of a non-polygonal domain that we experimented with was a standard disk. We discovered that the result of dropping a single grain in the center converges to some infinite tree inscribed in the disk (see Fig. \ref{fig_disk}). This tropical analytic curve has a deep geometric meaning -- it is the tropical caustic of the disk \cite{mikhalkin2023wave}. In general, tropical caustics of convex domains always emerge as such limits, in terms of our approach to establishing this result, we approximate the domain by considering its tropical wave front (level set of the tropical series defining the caustic) at a small time, which is always a polygon with rational slope sides.  

\subsection{Two perturbation points on a standard square}
Another direction of generalization of the basic example is through increasing the number of perturbation points. In this subsection, we deal with two perturbation points on the standard square domain. Denote these points by $\p$ and $\q.$ Our strategy for relaxing $\state{3}+\delta_\p+\delta_\q$ is the following: relax first $\state{3}+\delta_\p;$ then remove a grain from $\p,$ add grain at $\q,$ relax and repeat interchanging the roles of $\p$ and $\q$. It might happen that there is no need to repeat the process (see Fig. \ref{fig_twopointswaves}), i.e. when a wave send from $\q$ on $(\state{3}+\delta_\p)^\circ$ doesn't reach $\p.$ Otherwise, the edge of $(\state{3}+\delta_\p)^\circ$ passing through $\p$ can be moved from $\p$ (as on Fig. \ref{fig_twopointswaves2} ), i.e. $$((\state{3}+\delta_\p)^\circ-\delta_\p+\delta_\q)^\circ(\p)=3$$ and returning the grain at $\p$ results in an unstable state, therefore more iterations are needed.

\begin{figure}
\begin{tikzpicture}[scale=.55, every node/.style={scale=.8}]

\draw[step=0.5, very thin, gray!20](0,0) grid (10,10);

\draw(0,0)--(10,0)--(10,10)--(0,10)--(0,0);

\node at (5,1.5) {4};

\node at (0.5,0.5) {1};
\node at (1,1) {1};
\node at (1.5,1.5) {1};
\node at (2,2) {1};
\node at (2.5,2.5) {1};
\node at (3,3) {1};

\node at (0.5,9.5) {1};
\node at (1,9) {1};
\node at (1.5,8.5) {1};
\node at (2,8) {1};
\node at (2.5,7.5) {1};
\node at (3,7) {1};

\node at (9.5,0.5) {1};
\node at (9,1) {1};
\node at (8.5,1.5) {1};
\node at (8,2) {1};
\node at (7.5,2.5) {1};
\node at (7,3) {1};

\node at (9.5,9.5) {1};
\node at (9,9) {1};
\node at (8.5,8.5) {1};
\node at (8,8) {1};
\node at (7.5,7.5) {1};
\node at (7,7) {1};

\node at (3,3.5) {2};

\node at (3,4.5) {2};
\node at (3,5) {2};
\node at (3,5.5) {2};
\node at (3,6) {2};
\node at (3,6.5) {2};

\node at (7,3.5) {2};
\node at (7,4) {2};
\node at (7,4.5) {2};
\node at (7,5) {2};
\node at (7,5.5) {2};
\node at (7,6) {2};
\node at (7,6.5) {2};

\node at (3.5,3) {2};
\node at (4,3) {2};
\node at (4.5,3) {2};
\node at (5,3) {2};
\node at (5.5,3) {2};
\node at (6,3) {2};
\node at (6.5,3) {2};

\node at (3.5,7) {2};
\node at (4,7) {2};
\node at (4.5,7) {2};
\node at (5,7) {2};
\node at (5.5,7) {2};
\node at (6,7) {2};
\node at (6.5,7) {2};

\begin{scope}[xshift=330]
\draw[step=0.5, very thin, gray!20](0,0) grid (10,10);

\draw(0,0)--(10,0)--(10,10)--(0,10)--(0,0);

\node at (5,1.5) {4};

\node at (0.5,0.5) {2};
\node at (1,0.5) {0};

\node at (9.5,0.5) {2};
\node at (9,0.5) {0};

\node at (1.5,0.5) {2};
\node at (2,0.5) {2};
\node at (2.5,0.5) {2};
\node at (3,0.5) {2};
\node at (3.5,0.5) {2};
\node at (4,0.5) {2};
\node at (4.5,0.5) {2};
\node at (5,0.5) {2};
\node at (5.5,0.5) {2};
\node at (6,0.5) {2};
\node at (6.5,0.5) {2};
\node at (7,0.5) {2};
\node at (7.5,0.5) {2};
\node at (8,0.5) {2};
\node at (8.5,0.5) {2};

\node at (1.5,1) {1};
\node at (2,1.5) {1};
\node at (2.5,2) {1};
\node at (3,2.5) {1};

\node at (0.5,9.5) {1};
\node at (1,9) {1};
\node at (1.5,8.5) {1};
\node at (2,8) {1};
\node at (2.5,7.5) {1};
\node at (3,7) {1};

\node at (8.5,1) {1};
\node at (8,1.5) {1};
\node at (7.5,2) {1};
\node at (7,2.5) {1};

\node at (9.5,9.5) {1};
\node at (9,9) {1};
\node at (8.5,8.5) {1};
\node at (8,8) {1};
\node at (7.5,7.5) {1};
\node at (7,7) {1};

\node at (3,3) {2};
\node at (3,3.5) {2};

\node at (3,4.5) {2};
\node at (3,5) {2};
\node at (3,5.5) {2};
\node at (3,6) {2};
\node at (3,6.5) {2};

\node at (7,3) {2};
\node at (7,3.5) {2};
\node at (7,4) {2};
\node at (7,4.5) {2};
\node at (7,5) {2};
\node at (7,5.5) {2};
\node at (7,6) {2};
\node at (7,6.5) {2};

\node at (3.5,2.5) {2};
\node at (4,2.5) {2};
\node at (4.5,2.5) {2};
\node at (5,2.5) {2};
\node at (5.5,2.5) {2};
\node at (6,2.5) {2};
\node at (6.5,2.5) {2};

\node at (3.5,7) {2};
\node at (4,7) {2};
\node at (4.5,7) {2};
\node at (5,7) {2};
\node at (5.5,7) {2};
\node at (6,7) {2};
\node at (6.5,7) {2};
\end{scope}

\begin{scope}[yshift=-330]
\draw[step=0.5, very thin, gray!20](0,0) grid (10,10);

\draw(0,0)--(10,0)--(10,10)--(0,10)--(0,0);

\node at (5,1.5) {4};

\node at (0.5,0.5) {2};
\node at (1,0.5) {0};
\node at (1.5,0.5) {2};
\node at (1.5,1) {2};
\node at (2,1) {0};

\node at (9.5,0.5) {2};
\node at (9,0.5) {0};
\node at (8.5,0.5) {2};
\node at (8.5,1) {2};
\node at (8,1) {0};

\node at (2.5,1) {2};
\node at (3,1) {2};
\node at (3.5,1) {2};
\node at (4,1) {2};
\node at (4.5,1) {2};
\node at (5,1) {2};
\node at (5.5,1) {2};
\node at (6,1) {2};
\node at (6.5,1) {2};
\node at (7,1) {2};
\node at (7.5,1) {2};

\node at (2.5,1.5) {1};
\node at (3,2) {1};

\node at (0.5,9.5) {1};
\node at (1,9) {1};
\node at (1.5,8.5) {1};
\node at (2,8) {1};
\node at (2.5,7.5) {1};
\node at (3,7) {1};

\node at (7.5,1.5) {1};
\node at (7,2) {1};

\node at (9.5,9.5) {1};
\node at (9,9) {1};
\node at (8.5,8.5) {1};
\node at (8,8) {1};
\node at (7.5,7.5) {1};
\node at (7,7) {1};

\node at (3,2.5) {2};
\node at (3,3) {2};
\node at (3,3.5) {2};

\node at (3,4.5) {2};
\node at (3,5) {2};
\node at (3,5.5) {2};
\node at (3,6) {2};
\node at (3,6.5) {2};

\node at (7,2.5) {2};
\node at (7,3) {2};
\node at (7,3.5) {2};
\node at (7,4) {2};
\node at (7,4.5) {2};
\node at (7,5) {2};
\node at (7,5.5) {2};
\node at (7,6) {2};
\node at (7,6.5) {2};

\node at (3.5,2) {2};
\node at (4,2) {2};
\node at (4.5,2) {2};
\node at (5,2) {2};
\node at (5.5,2) {2};
\node at (6,2) {2};
\node at (6.5,2) {2};

\node at (3.5,7) {2};
\node at (4,7) {2};
\node at (4.5,7) {2};
\node at (5,7) {2};
\node at (5.5,7) {2};
\node at (6,7) {2};
\node at (6.5,7) {2};

\end{scope}

\begin{scope}[yshift=-330,xshift=330]
\draw[step=0.5, very thin, gray!20](0,0) grid (10,10);

\draw(0,0)--(10,0)--(10,10)--(0,10)--(0,0);

\node at (5,1.5) {2};

\node at (0.5,0.5) {2};
\node at (1,0.5) {0};
\node at (1.5,0.5) {2};
\node at (1.5,1) {2};
\node at (2,1) {0};
\node at (2.5,1) {2};
\node at (3,1.5) {0};

\node at (9.5,0.5) {2};
\node at (9,0.5) {0};
\node at (8.5,0.5) {2};
\node at (8.5,1) {2};
\node at (8,1) {0};
\node at (7.5,1) {2};
\node at (7,1.5) {0};

\node at (2.5,1.5) {2};
\node at (3,2) {2};

\node at (0.5,9.5) {1};
\node at (1,9) {1};
\node at (1.5,8.5) {1};
\node at (2,8) {1};
\node at (2.5,7.5) {1};
\node at (3,7) {1};

\node at (7.5,1.5) {2};
\node at (7,2) {2};

\node at (9.5,9.5) {1};
\node at (9,9) {1};
\node at (8.5,8.5) {1};
\node at (8,8) {1};
\node at (7.5,7.5) {1};
\node at (7,7) {1};

\node at (3,2.5) {2};
\node at (3,3) {2};
\node at (3,3.5) {2};

\node at (3,4.5) {2};
\node at (3,5) {2};
\node at (3,5.5) {2};
\node at (3,6) {2};
\node at (3,6.5) {2};

\node at (7,2.5) {2};
\node at (7,3) {2};
\node at (7,3.5) {2};
\node at (7,4) {2};
\node at (7,4.5) {2};
\node at (7,5) {2};
\node at (7,5.5) {2};
\node at (7,6) {2};
\node at (7,6.5) {2};

\node at (3.5,1.5) {1};
\node at (4,1.5) {1};
\node at (4.5,1.5) {1};

\node at (5.5,1.5) {1};
\node at (6,1.5) {1};
\node at (6.5,1.5) {1};

\node at (3.5,7) {2};
\node at (4,7) {2};
\node at (4.5,7) {2};
\node at (5,7) {2};
\node at (5.5,7) {2};
\node at (6,7) {2};
\node at (6.5,7) {2};

\end{scope}

\end{tikzpicture}
\caption{A relaxation seen as a decomposition into three waves where we start with the middle picture of Fig. \ref {fig_relaxone34} and add another grain at some point (marked with value $4$ on the first three pictures). Each wave deforms the tropical curve so that the face containing the new perturbation point shrinks towards it, collapsing to an edge after a third wave. This edge is of multiplicity two since it is horizontal but made of values $1$ and not $2,$ i.e. the deviation from $3$ is doubled. Such an edge is not made of a soliton pattern, but a superposition of two, i.e. a multi-soliton.}
\label{fig_twopointswaves}
\end{figure}

The reader might wonder, why we prefer to perform the relaxation this way, by removing a grain at one point and putting it back at the other. The answer is that this way the deviation locus is made only of solitons and triads (if we are not unlucky, this subtlety will be discussed in the next paragraph), having an extra grain on one such edge would count as a defect which would let the waves through it. In other words, by removing grains at all but one perturbation point $\p$ we ensure that the waves are confined to the face (in the complement of a tropical curve) containing $\p$ and act simply by shrinking this face until its boundary passes through $\p$. At the tropical level, we call this operator $G_\p.$

\begin{figure}
\begin{tikzpicture}[scale=.55, every node/.style={scale=.8}]
\draw[step=0.5, very thin, gray!20](0,0) grid (10,10);

\draw(0,0)--(10,0)--(10,10)--(0,10)--(0,0);

\node at (4,4.5) {4};

\node at (0.5,0.5) {1};
\node at (1,1) {1};
\node at (1.5,1.5) {1};
\node at (2,2) {1};
\node at (2.5,2.5) {1};
\node at (3,3) {1};

\node at (0.5,9.5) {1};
\node at (1,9) {1};
\node at (1.5,8.5) {1};
\node at (2,8) {1};
\node at (2.5,7.5) {1};
\node at (3,7) {1};

\node at (9.5,0.5) {1};
\node at (9,1) {1};
\node at (8.5,1.5) {1};
\node at (8,2) {1};
\node at (7.5,2.5) {1};
\node at (7,3) {1};

\node at (9.5,9.5) {1};
\node at (9,9) {1};
\node at (8.5,8.5) {1};
\node at (8,8) {1};
\node at (7.5,7.5) {1};
\node at (7,7) {1};

\node at (3,3.5) {2};

\node at (3,4.5) {2};
\node at (3,5) {2};
\node at (3,5.5) {2};
\node at (3,6) {2};
\node at (3,6.5) {2};

\node at (7,3.5) {2};
\node at (7,4) {2};
\node at (7,4.5) {2};
\node at (7,5) {2};
\node at (7,5.5) {2};
\node at (7,6) {2};
\node at (7,6.5) {2};

\node at (3.5,3) {2};
\node at (4,3) {2};
\node at (4.5,3) {2};
\node at (5,3) {2};
\node at (5.5,3) {2};
\node at (6,3) {2};
\node at (6.5,3) {2};

\node at (3.5,7) {2};
\node at (4,7) {2};
\node at (4.5,7) {2};
\node at (5,7) {2};
\node at (5.5,7) {2};
\node at (6,7) {2};
\node at (6.5,7) {2};

\begin{scope}[xshift=330]
\draw[step=0.5, very thin, gray!20](0,0) grid (10,10);

\draw(0,0)--(10,0)--(10,10)--(0,10)--(0,0);

\node at (3,4) {4};
\node at (4,4.5) {4};

\node at (0.5,0.5) {1};
\node at (1,1) {1};
\node at (1.5,1.5) {1};
\node at (2,2) {1};
\node at (2.5,2.5) {1};
\node at (3,3) {1};
\node at (3.5,3.5) {1};
\node at (4,4) {1};

\node at (0.5,9.5) {1};
\node at (1,9) {1};
\node at (1.5,8.5) {1};
\node at (2,8) {1};
\node at (2.5,7.5) {1};
\node at (3,7) {1};
\node at (3.5,6.5) {1};
\node at (4,6) {1};

\node at (9.5,0.5) {1};
\node at (9,1) {1};
\node at (8.5,1.5) {1};
\node at (8,2) {1};
\node at (7.5,2.5) {1};
\node at (7,3) {1};
\node at (6.5,3.5) {1};
\node at (6,4) {1};

\node at (9.5,9.5) {1};
\node at (9,9) {1};
\node at (8.5,8.5) {1};
\node at (8,8) {1};
\node at (7.5,7.5) {1};
\node at (7,7) {1};
\node at (6.5,6.5) {1};
\node at (6,6) {1};

\node at (4,5) {2};
\node at (4,5.5) {2};

\node at (6,4.5) {2};
\node at (6,5) {2};
\node at (6,5.5) {2};

\node at (4.5,4) {2};
\node at (5,4) {2};
\node at (5.5,4) {2};

\node at (4.5,6) {2};
\node at (5,6) {2};
\node at (5.5,6) {2};

\end{scope}

\begin{scope}[yshift=-330]
\draw[step=0.5, very thin, gray!20](0,0) grid (10,10);

\draw(0,0)--(10,0)--(10,10)--(0,10)--(0,0);

\node at (3,4) {2};
\node at (4,4.5) {4};

\node at (1,2.5) {2};
\node at (1,3) {2};
\node at (1,3.5) {2};
\node at (1,4) {2};
\node at (1,4.5) {2};
\node at (1,5) {2};
\node at (1,5.5) {2};
\node at (1,6) {2};
\node at (1,6.5) {2};
\node at (1,7) {2};
\node at (1,7.5) {2};

\node at (0.5,0.5){2};
\node at (0.5,1){0};
\node at (0.5,1.5){2};
\node at (1,1.5) {2};
\node at (1,2) {0};

\node at (0.5,9.5){2};
\node at (0.5,9){0};
\node at (0.5,8.5){2};
\node at (1,8.5) {2};
\node at (1,8) {0};

\node at (1.5,2.5) {1};
\node at (2,3) {1};
\node at (2.5,3.5) {1};

\node at (1.5,7.5) {1};
\node at (2,7) {1};
\node at (2.5,6.5) {1};
\node at (3,6) {1};

\node at (9.5,0.5) {1};
\node at (9,1) {1};
\node at (8.5,1.5) {1};
\node at (8,2) {1};
\node at (7.5,2.5) {1};
\node at (7,3) {1};
\node at (6.5,3.5) {1};
\node at (6,4) {1};

\node at (9.5,9.5) {1};
\node at (9,9) {1};
\node at (8.5,8.5) {1};
\node at (8,8) {1};
\node at (7.5,7.5) {1};
\node at (7,7) {1};
\node at (6.5,6.5) {1};
\node at (6,6) {1};

\node at (3,4.5) {2};
\node at (3,5) {2};
\node at (3,5.5) {2};

\node at (6,4.5) {2};
\node at (6,5) {2};
\node at (6,5.5) {2};

\node at (3.5,4) {2};
\node at (4,4) {2};
\node at (4.5,4) {2};
\node at (5,4) {2};
\node at (5.5,4) {2};

\node at (3.5,6) {2};
\node at (4,6) {2};
\node at (4.5,6) {2};
\node at (5,6) {2};
\node at (5.5,6) {2};
\end{scope}

\begin{scope}[yshift=-330,xshift=330]
\draw[step=0.5, very thin, gray!20](0,0) grid (10,10);

\draw(0,0)--(10,0)--(10,10)--(0,10)--(0,0);

\node at (3,4) {2};

\node at (1,2.5) {2};
\node at (1,3) {2};
\node at (1,3.5) {2};
\node at (1,4) {2};
\node at (1,4.5) {2};
\node at (1,5) {2};
\node at (1,5.5) {2};
\node at (1,6) {2};
\node at (1,6.5) {2};
\node at (1,7) {2};
\node at (1,7.5) {2};

\node at (0.5,0.5){2};
\node at (0.5,1){0};
\node at (0.5,1.5){2};
\node at (1,1.5) {2};
\node at (1,2) {0};

\node at (0.5,9.5){2};
\node at (0.5,9){0};
\node at (0.5,8.5){2};
\node at (1,8.5) {2};
\node at (1,8) {0};

\node at (1.5,2.5) {1};
\node at (2,3) {1};
\node at (2.5,3.5) {1};
 
\node at (3.5,4.5) {1};

\node at (1.5,7.5) {1};
\node at (2,7) {1};
\node at (2.5,6.5) {1};
\node at (3,6) {1};
\node at (3.5,5.5) {1};

\node at (9.5,0.5) {1};
\node at (9,1) {1};
\node at (8.5,1.5) {1};
\node at (8,2) {1};
\node at (7.5,2.5) {1};
\node at (7,3) {1};
\node at (6.5,3.5) {1};
\node at (6,4) {1};
\node at (5.5,4.5) {1};

\node at (9.5,9.5) {1};
\node at (9,9) {1};
\node at (8.5,8.5) {1};
\node at (8,8) {1};
\node at (7.5,7.5) {1};
\node at (7,7) {1};
\node at (6.5,6.5) {1};
\node at (6,6) {1};
\node at (5.5,5.5) {1};

\node at (3.5,5) {2};

\node at (5.5,5) {2};

\node at (4.5,4.5) {2};
\node at (5,4.5) {2};

\node at (4,5.5) {2};
\node at (4.5,5.5) {2};
\node at (5,5.5) {2};

\end{scope}

\end{tikzpicture}
\caption{Again, start with the middle picture of Fig. \ref {fig_relaxone34} and add another grain at $\q$ two steps to the right and one step up from the original perturbation point $\p$. The second picture is derived from the first by removing a grain from $\p,$ sending two waves from $\q$ and adding back a grain at $\p$. Next, we remove a grain from $\q,$ send two waves from $\p$ and add the grain at $\q.$ The last picture, the result of the relaxation for $\langle 3\rangle+\delta_{\p}+\delta_{\q},$ is obtained by sending a single wave from $\q.$   }
\label{fig_twopointswaves2}
\end{figure}

Now, about that subtlety of being unlucky in the previous paragraph. Ideally, we want to have perfect control of all the relaxations. Since we decompose our relaxations into waves, this is guaranteed if deviation loci are always made only of solitons and triads (because we know how they move under wave action). As the example of Fig. \ref{fig_twopointswaves} shows, this is a utopia: sometimes the results of relaxations for states of our interest are made of more complicated patterns that include multi-solitons or higher than nodal singularities (the theory of which is not yet developed). However, such a problem might only arise at terminal stages of sending waves from a given perturbation point, when the corresponding face of a tropical curve shrinks to an edge or a vertex, and therefore it can be avoided when proving the main theorem, i.e. by refraining from sending those last few waves at each step, we still have a lower bound on the toppling function.

Another remark is that the way we approached this example of two points on the square is not exactly how we would treat it in the proof of the scaling limit theorem. The main difference lies in the idea that we don't want our waves to interact with junctions of the boundary of the domain and solitons. Although we know what to expect from such interactions, and for small slopes of solitons and boundary sides it can be verified by hand, a general rigorous theory for such situations is missing. To circumvent that, similar to the above case of a single point perturbation on a unimodular polygon, we replace the maximal stable state with some auxiliary state $\psi$ with deviation locus concentrated along the boundary. Here again, it is made of solitons joining in triads and going along the sides of the square, but to guarantee that no wave in the present example reaches the boundary we actually need two solitons going along one of the sides. To decide which side is that we look first at the expected limit of the toppling function, which is a tropical polynomial $f_{P;\Omega}$ vanishing on the boundary of the square $\Omega.$ An important characteristic of such an  $\Omega$-tropical polynomial is its quasi-degree which measures its slope along each side against the slope of the primitive support function of this side. For a one-point perturbation on a unimodular polygon, the quasi-degree takes the value one for each of its sides -- this is why the state $\psi$ had a single soliton along each side in that case. In the case of a two-point perturbation, the quasi-degree of $f_{P;\Omega}$ may take the value two for one side of the polygon, therefore we need two solitons going along this side to guarantee that waves in the relaxation of $\psi+\delta_p+\delta_q$ never reach the boundary and interact only with solitons, triads and nodes.

\subsection{Summary of the strategy} The generality of the present article requires combining all the ideas encountered above and a few more tricks. For an admissible domain $\Omega$ and a finite collection of points $P\subset\Omega^\circ,$ and $\h>0$ we consider a perturbation $\phi_\h$ of the maximal stable state (all-three-state) on $\Omega_\h=\Omega^\circ\cap \h\mathbb{Z}^2$ by adding a single grain at roundings of points in $P.$ We give upper and lower bounds on the toppling function of $\phi_\h$ multiplied by $\h$ that both converge to the $\Omega$-tropical series $f_{\Omega,P}=G_P 0_\Omega,$ which is the minimal among those not linear at points $P.$ The upper bound is the easier one and is achieved through the least action principle in sandpiles, but here also lies the origin of a technical condition of properness for roundings -- in essence, we require the roundings to belong to the non-harmonicity locus of the discretization of $f_{\Omega,P}.$ Combining the two bounds, it is not hard to deduce the (Haussdorf, on compacts) convergence of the deviation loci $D(\phi^\circ)$ to the tropical curve $C(f_{\Omega,P}),$ this is done in Section \ref{sec_proof}.

The lower bound is more involved and is achieved through a variety of approximations of diverse origins. First, we reduce, in an almost canonical way, the case of a general admissible domain $\Omega$ to a $\mathbb{Q}$-polygon contained in it so that the expected limits of the toppling functions differ by a small constant (Lemma \ref{lemma_epsfun}). Such a $\mathbb{Q}$-polygon and the expected limit of the toppling function on it still may have unwanted pathologies, so we further cut out some corners of the polygon (``blow-up'' in toric terminology) until it becomes unimodular and the limiting tropical polynomial has a nice quasi-degree (Definitions \ref{def_quasidegree} and \ref{def_nice}). These conditions, by Proposition \ref{prop_smoothnice}, guarantee the existence of a tropical polynomial $g$ with the same quasi-degree, such that its tropical curve $C(g)$ is smooth and belongs to a small neighborhood of the boundary. We use this polynomial $g$ to carve an auxiliary state $\psi_\h$ made of solitons and triads (see Section \ref{sec_smoothings}) going along $C(g)$ by which we substitute the maximal stable state. Since $G_P g$ has the same quasi-degree as $g,$ no wave will interact with the boundary in the corresponding relaxation.

Working with the perturbation of this auxiliary state may be thought of as a partial relaxation of the original $\phi_h:$ on the one hand, we don't perform topplings outside the domain of definition of $\psi_\h;$ and on the other hand, we do less topplings along its deviation locus. We further reduce this lower bound due to the following two technicalities. In the continuous version of the sandpile model, the passage from $g$ to $G_P g$ might require an infinite sequence of applications of the single grain operators $G_q,$ for $q\in P,$ which, however, converges to $G_P g.$ For our purposes, we don't need to know if indeed such infinite sequences exist, or all of them stabilize, and we simply choose a particular $G=G_{q_m}G_{q_{m-1}}\dots G_{q_1}$ with $q_j \in P$ such that $0\leq G_Pg-Gg\leq \varepsilon$ and that each application of $G_{q_j}$ in $Gg$ is effective, i.e. every operator modifies the tropical polynomial (otherwise it can be omitted from $G$). At the discrete level, we will follow this sequence $q_1,\dots, q_m$ by ``activating'' points one by one and sending some number of waves from the corresponding roundings. This is the first technicality reducing the lower bound, i.e. we don't follow the whole, potentially infinite, tropical relaxation at the discrete level.

The second technicality is that we send less waves at each individual round of activating a point $q_j$ so that we can stay within the realm of solitons, triads and nodes. Explicitly, being effective in $g\mapsto Gg$ for an operator $G_{q_j}$ means that it acts by increasing a $v_j$-coefficient by some $c_j>0,$ i.e. it acts by $\Add_{v_j}^{c_j}.$ Now, we replace each $G_{q_j}$ in $G$ by $\Add_{v_j}^{c_j-\eta}$ for some $\eta>0.$ This way, a face of the tropical curve is never shrank and we encounter only at most nodal curves while applying $\Add_{v_j}^{t(c_j-\eta)},$ for $t\in[0,1]$, by Proposition \ref{prop_atmostnodal}. Here we use our perfect control paradigm, where each individual wave, by interaction with solitons, simply increases a coefficient of the corresponding tropical polynomial, Lemma \ref{lemma_wavegp}. Through this process, the partial toppling function multiplied by $\h$ will be at least $Gg-g-m\eta,$ so we may take $\eta$ to be $m^{-1}\varepsilon$ (recall that $m$ is the number of operators $G_{q_{j}}$ in $G$).

\section{Proper roundings and the upper bound}
\label{sec_proproundestimabove}
Recall that $\Omega$ denotes an admissible convex domain, $P$ is a finite subset in its interior and $f_{\Omega,P}=G_P0_\Omega$ is the minimal $\Omega$-tropical series such that it is not smooth at all points of $P.$ This series is an expected limit of the toppling function for the perturbation of the maximal stable state at proper roundings of points from $P.$

\begin{definition}
\label{def_proper}
For $\h>0,$ a set of roundings $P^\h = \{\p^\h|\p\in P\}$ for the set of points $P$ is called {\it proper} if the function 
\begin{equation}
\label{eq_estimate}
F:\Omega_\h\to\ZZ_{\geq 0}, F(z)=[\h^{-1}f_{\Omega,P}(z)]
\end{equation}
 has negative discrete Laplacian at all points $\p^\h$. Here $[\cdot]$
stands for the usual integer part of a non-negative number.
\end{definition}

\begin{proposition}
\label{prop_roundings}
For each finite subset $P$ of $\Omega^\circ$ there exists a set $P^\h$ of proper roundings. 
\end{proposition}

\begin{proof}
\label{proof_proproundings} 
Recall that the piece-wise linear function $f_{\Omega,P}:\Omega\to\RR_{\geq 0}$ is not smooth at all points in $P$. Let $$f_{\Omega,P}(x,y)=\inf_{(i,j)\in\A} (ix+jy+a_{ij}).$$ We consider the function \eqref{eq_estimate}. Note that on $\Omega_\h$ we have
\begin{equation}
\label{eq_rounded}
\h F(x,y) = \min_{(i,j)\in\A} (ix+jy+\h[a_{ij}\h^{-1}]). 
\end{equation}

The difference between corresponding coefficients $a_{ij}$ and $\h[a_{ij}\h^{-1}]$ is at most $\h$. 
It follows from Lemma~\ref{lemma_ecloseseries} below and Remark after it that for each $\p\in C(f_{\Omega,P})$ there exists an $\h$-close to $\p$ point $\p^\h\in\Omega_\h$ such that $\p^\h$ and one of its neighbors in $\Omega_\h$ belong to different regions of linearity of $F$. This implies that $\Delta F(\p^\h_i)<0$ for $i=1,\dots n$.
\end{proof}

\begin{lemma}[\cite{us_series}]
\label{lemma_ecloseseries}
Let $\e>0,\A\subset\ZZ^2$ and $f,g$ be two $\Omega$-tropical series  written as
$$f(x,y)=\inf_{(i,j)\in\A}(ix+jy+a_{ij}),g(x,y)=\inf_{(i,j)\in\A}(ix+jy+a_{ij}+\eta_{ij}).$$
If $|\eta_{ij}|<\e$  for each $(i,j)\in\A$, then $C(f)$ is $2\e$-close to $C(g)$.
\end{lemma}
\begin{remark}
A proof of this lemma goes as follows: changing $a_{ij}$ to $a_{ij}+\eta_{ij}$ moves edges of our $\Omega$-tropical series by at most $\eta_{ij}$. If all the numbers $\eta_{ij}$ are of the same sign, then each edge of the tropical curve moves twice, but in opposite directions. Thus, we will use the updated version of this lemma: if all $\eta_{ij}$ are of the same sign and $|\eta_{ij}|<\e$ for each $(i,j)\in\A$, then $C(f)$ is $\e$-close to $C(g)$.   
\end{remark}

In fact, the set $P^\h=\{\p^\h\}$ of proper roundings depends on $\Omega$ and $P$, so we should write $\p^\h_{\Omega,P}$ for each point $\p$. 
Nevertheless, for a fixed $\h$ small enough this
rounding $\p^\h$ of $\p\in P$ depends only on the behavior of $C(f_{\Omega,P})$ in a small
neighborhood of $\p$. The choice in the above proof depends only on an arbitrary small neighborhood of $\p\in P$ on $C(\Omega,P)$. Therefore we can fix choices (for example: ``take the nearest point in $\Omega_\h$ from the south-east region of $\p$'') for all possible neighbors of a points in a tropical curve.

\begin{corollary}
\label{cor_roundings2}
If $\p\in P\cap P', \p\in \Omega\cap\Omega'$ and $C(f_{\Omega,P})$ coincides
with $C(f_{\Omega',P'})$ in a neighborhood of $\p$, then we can take the corresponding proper roundings to be equal,  
$\p^\h_{\Omega,P}=\p^\h_{\Omega',P'}$.   
\end{corollary}


\begin{proposition}
\label{prop_abovebound}
The function $$F(z)=[\h^{-1}f_{\Omega,P}(z)], z\in\Omega_\h$$ bounds from
above the toppling function $H_{\phi_\h}$ of $\phi_\h=\state{3}+\sum_{\p\in P}\delta_{\p^\h}$. 
\end{proposition}
\begin{proof}
\label{proof_propabovebound}
We apply the Least Action Principle (see \cite{FLP,sadhu2011effect} for finite graphs, \cite{us_solitons} for infinite graphs and locally finite relaxations): if $F\geq 0$ on a graph and $$\phi+\Delta F\leq \state{3},$$ then $F$ bounds the toppling function of $\phi$ from above. In our case, proper roundings exist and $\Delta F(\p^\h)<0$ for $\p\in P$ by the definition.
\end{proof}

On the other hand, in a certain setting a rounding is not necessary at all, as it was in \cite{announce}, Section 1.2, where $\Omega$ is a lattice polygon and $P$ belongs to the lattice.
\begin{proposition}
\label{prop_roundinglattice}
If $P\subset \ZZ^2$, $\Omega$ is lattice polygon, and $\h^{-1}\in\NN$, then we can take the proper roundings $\p^\h=\p$ for each $\p\in P$.
\end{proposition}

\begin{proof}
\label{proof_roundinglattice}
If $P\subset \ZZ^2,\h^{-1}\in\NN$ and $\Omega$ is a lattice polygon, then for $f=f_{\Omega,P}$, in \eqref{eq_series} $a_{ij}\in \ZZ$. Indeed, near the boundary of $\Omega$ that holds because $\Omega$ is a lattice polygon, and then when a linear function $ix+jy+a_{ij}$ with $a_{ij},i,j\in\ZZ$ is equal to another linear function $i'x+j'y+a_{i'j'}, i',j'\in \ZZ^2$ at $\p\in \ZZ^2$ this guarantees that $a_{i'j'}$ is also integer. Therefore $\h[\h^{-1}a_{ij}]=a_{ij}$ in the proof of the Proposition~\ref{prop_roundings}, so we can take $\p^\h=\p$ for all $\p\in P$. 
\end{proof}



\section{Sandpiles on $\QQ$-polygons}

\begin{definition}
Let $\Pi\subset\RR^2$ be a finite intersection of a non-empty collection of half-planes with rational slopes.
We call $\Pi$ a {\it $\QQ$-polygon} if it is a closed set with a non-empty interior.
\end{definition}

If $\Omega=\Pi$ is a $\QQ$-polygon, then we have good control over the behavior of $D(\phi^\circ_\h)$ by means of \cite{us_solitons}, and reducing the general problem to $\QQ$-polygons seems to be an unavoidable technical step. 
Note that a $\QQ$-polygon is not necessarily compact. It is easy to verify that a $\QQ$-polygon is admissible (Definition~\ref{def_omegaadmissible}). The next lemma provides us with a family of $\QQ$-polygons approximating $\Omega$ from inside that we use to extend the proof of the main result from $\QQ$-polygons to general admissible domains.
\begin{lemma}[\cite{us_series}]
\label{lemma_epsfun} For each compact set $K\subset \Omega^\circ$ such that $P\subset K$ and for each $\e>0$ small enough there exists a $\QQ$-polygon $\Omega_{\e,K}\subset \Omega$ such that $B_{3\e}(K)\subset \Omega_{\e,K}$ and the following holds: 
$$f_{\Omega,P}= f_{\Omega_{\e,K},P}+\e \text{\ on\ } B_{3\e}(K).$$ 
\end{lemma}
Here $B_{3\e}(K)$ denotes the $3\e$-neighbourhood of $K.$ The idea of a proof: if $\Omega$ is a compact set, then $\Omega_{\e,K}$ is the set $\{f_{\Omega,P}\geq \e\}$ for $\e$ small enough; for non-compact admissible $\Omega$ we additionally cut the set $\{f_{\Omega,P}\geq \e\}$ far enough from $K$.

\begin{corollary}
\label{cor_epscurve} Tropical curves defined by $f_{\Omega,P} $ and $f_{\Omega_{\e,K},P}$ coincide on $K$,
i.e.  $$C(f_{\Omega,P})\cap K=C(f_{\Omega_{\e,K},P})\cap K.$$
\end{corollary}

 \begin{proposition}
\label{prop_qpolygonconvergence}
Suppose that $\Omega=\Pi\subset \RR^2$ is a $\QQ$-polygon. Choose any $\e>0$. Then, for all $\h>0$ small enough, the toppling functions $H_{\phi_h}$
of the states $\phi_\h=\state{3}+\sum_{\p\in P}\delta_{\p^\h}$  satisfy
$$\h H_{\phi_h}(z)
>f_{\Pi,P}(z)-\e \text{\ at\ all\ } z\in \Pi_\h.$$ 
\end{proposition}
This Proposition is a cornerstone of all our results, we prove it in Section~\ref{sec_lower} by assembling (from so-called sandpile solitons and triads reviewed in Section \ref{sec_smoothings}) an auxiliary state $\psi_\h$ discretizing a smooth tropical curve near the boundary (see Section \ref{sec_series}) with the property $$\h H_{\phi_h}\geq \h H_{\psi_h}> f_{\Pi,P}-\e.$$


\section{Proof of the main theorem}
\label{sec_proof}


%
%

\begin{proof}[Proof of Theorem~\ref{th_main}]
\label{proof_maintheorem}
Consider a compact set $K\subset\Omega^\circ$, such that $P\subset K$, and choose any $\e>0$. We are going to prove that $D(\phi_\h^\circ)\cap K$ is $3\e$-close to $C(f_{\Omega,P})\cap K$ as long as $\h$ is small enough.

Choose a $\QQ$-polygon $\Pi^\e=\Omega_{\e,K}$ by Lemma~\ref{lemma_epsfun}. Consider the state $$\phi_\h^\e=\langle 3 \rangle+\sum_{\p\in P}
\delta_{\p^\h}$$ on $\Pi^\e_\h = \Pi^\e\cap\h\ZZ^2$. 

Note that the proper roundings $\p^\h$ of points $\p\in P$
on $\Pi^\e_\h$ can be taken the same as on $\Omega_\h$ by Corollary
\ref{cor_epscurve} and Corollary~\ref{cor_roundings2}. 

Denote by $H$ the toppling function of $\phi_\h$ (on $\Omega_\h$) and by $H^\e$ the toppling function of $\phi_\h^\e$ (on $\Pi^\e_\h$). Note that $(\phi_\h^\e)^\circ$ can be thought of as a partial relaxation of $\phi_\h,$ wher we don't do topplings at  vertices outside $\Pi^\e_\h$; 
therefore, $\h H\geq \h H^\e$, (2) below. 

Since $\Pi^\e$ is a $\QQ$-polygon, then, by Proposition~\ref{prop_qpolygonconvergence} we can
choose $\h$ small enough, such that $\h H^\e>f_{\Pi^\e,P}-\e$, (1) below.

On $B_{3\e}(K)$, combining the above arguments with Proposition~\ref{prop_abovebound}, (3) below, and Lemma~\ref{lemma_epsfun}, (5) below, we obtain
 \begin{equation}
\label{eq_estimatetoppling}
f_{\Pi^\e,P}-\e\stackrel{(1)}{<}\h H^\e\stackrel{(2)}{\leq} \h H
\stackrel{(3)}{\leq} \h[\h^{-1}f_{\Omega,P}]\stackrel{(4)}{\leq} f_{\Omega,P}\stackrel{(5)}{=} f_{\Pi^\e,P}+\e. 
\end{equation}

Hence $|H-H^\e|\leq\frac{2\e}{\h}$ and by Lemma~\ref{lemma_deviation} below the deviation set $D(\phi_\h^\circ)$ is $2\e$-close to
$$(D\big((\phi_\h^\e)^\circ\big)\cap K)\cup \partial B_{3\e}(K),$$
on $B_{3\e}(K)$. 

\begin{lemma}
\label{lemma_deviation}
If the toppling function $H_\psi$ of a state $\psi$ on $\Omega_\h$ is bounded by a constant
$\C\in \ZZ_{>0}$, then $$D(\psi^\circ)\subset B_{\C\cdot\h}(D(\psi)\cup\partial\Omega),$$
where $B_r(X)$ is the $r$-neighborhood of a set $X\subset \RR^2$.
\end{lemma}

\begin{proof}
\label{proof_lemmadeviation}

Consider a point $z\in D(\psi^\circ)$. Suppose that $z$ does not belong to $D(\psi)$ or $\partial \Omega_\h$. Since $\psi^\circ(z)<\psi(z)$, we have $\Delta H_{\psi}(z)< 0$, therefore there exists a neighbor $z_1$ of $z$ such that $H_{\psi}(z_1)< H_{\psi}(z)$. If $z_1$ does not belong to $D(\psi)$ or $\partial\Omega_\h$, then $\Delta H_{\psi}(z_1)\leq 0$, and $H_{\psi}(z_1)< H_{\psi}(z)$ implies that $z_1$ has a neighbor $z_2$ such that $H_{\psi}(z_2)<H_{\psi}(z_1)$. We repeat this argument and find $z_3,z_4,$ etc. Since $0\leq H_{\psi}\leq \C$, we can not have such a chain of length bigger than $\C+1$. Therefore, starting with any point $z\in D(\psi^\circ)$ and passing each time to a neighbor we reach $D(\psi)$ or $\partial\Omega_\h$ by at most $\C$ steps, which concludes the proof. 
\end{proof}

By Proposition~\ref{prop_qpolygonconvergence}, $D((\phi_\h^\e)^\circ)\cap K$ is $\e$-close to
$C(f_{\Pi^\e,P})\cap K$ which is, in turn, coincides with $C(f_{\Omega,P})\cap K$ (Corollary~\ref{cor_epscurve}). Thus, we proved that  $D(\phi_\h^\circ)\cap K$ is $3\e$-close to $C(f_{\Omega,P})\cap K$.

It is left to explain that $D(\phi_\h^\circ)$ doesn't have holes with respect to $C(f_{\Omega,P}),$ i.e. every point of $C(f_{\Omega,P})$ is approximated by a point of $D(\phi_\h^\circ).$ Such a hole cannot exist since (a) the toppling function $H_{\phi_\h}$ is linear on the sufficiently big harmonicity regions, (b) its gradient inside each face of $C(f_{\Omega,P})$ is the same as the gradient of $f_{\Omega,P}$ on this face (the later gradients are different for different faces). In other words, having a global hole in $D(\phi_\h^\circ)$  would mean that one can pass between two different faces staying within the linearity region, thus, equating two different gradients on the two sides of the hole. To show (b) it is enough to observe that otherwise, the toppling function rescaled by $\h$ couldn't be between $f_{\Omega,P}$ and $f_{\Omega,P}-\e$ on a big enough region of linearity because the gradients are integral. To show (a), apply twice the following classical theorem of Duffin \cite{Duffin}:
\begin{theorem*}
Let $R>1,$ $v\in\mathbb{Z}^2,$ $\Gamma=B_R(v)\cap\mathbb{Z}^2$ and $f\colon\Gamma\rightarrow\mathbb{R}$ be a non-negative harmonic function. Let $v'\in\mathbb{Z}^2,$ such that $|v'-v|=1.$ Then,
$$|f(v')-f(v)|\leq \varkappa R^{-1} \max_{w\in\Gamma}f(w),$$ where $\varkappa$ is an absolute constant.
\end{theorem*} 
 
In more details, assume that $B_\gamma(v)\subset\Omega\backslash D(\phi_\h^\circ)$ for some $\gamma>0$ and $v\in\Omega_\h.$ Then, $$|\partial_\bullet H_{\phi_\h}(v)|\leq \varkappa(h^{-1}\gamma)^{-1}\h^{-1}\mu=\varkappa\gamma^{-1}\mu,$$ where $\mu$ is the maximum of $f_{\Omega,P}$ and $\partial_\bullet$ is some first discrete derivative (either with respect to horizontal or vertical direction). To apply the theorem twice for $w\in B_\gamma(v)$ with the later upper bound, we need to assume that $B_{2\gamma}(v)\subset\Omega\backslash D(\phi_\h^\circ)$ and make $\partial_\bullet H_{\phi_\h}$ non-negative by adding $\varkappa\gamma^{-1}\mu,$ that gives $$|\partial_\bullet\partial_\bullet H_{\phi_\h}(w)|=|\partial_\bullet(\partial_\bullet H_{\phi_\h}(w)+\varkappa\gamma^{-1}\mu)|\leq 2\varkappa^2\gamma^{-2}\mu\h.$$ If we keep the same $\gamma$ and take $h$ sufficiently small, this gives $$|\partial_\bullet\partial_\bullet H_{\phi_\h}(w)|<1$$ for all $w\in B_\gamma(v).$ Since the toppling function is integer-valued, we have the vanishing of all second derivatives on globally observable harmonicity regions.


\end{proof}

\begin{remark}
\label{rem_th5}
Note that \eqref{eq_estimatetoppling} implies that $f_{\Omega,P}=\lim_{\h\to 0} \h H_{\phi_\h}$. This gives a stronger version of Theorem~5 announced in \cite{announce}.
\end{remark}

\section{Tropical series and the dynamic by operators $G_\p$}\label{sec_series}

The following three sections are dedicated to a proof of Proposition~\ref{prop_qpolygonconvergence}. The theory of tropical series is developed in \cite{us_series}. Here we collect the results and definitions from \cite{us_series} that we need for our proofs.

Consider an admissible domain $\Omega$ and $P=\{\p_1,\dots,\p_n\}$, a
finite collection of points in $\Omega^\circ$. Let $f$ be an $\Omega$-tropical series. 

\begin{definition}[\cite{us_series}]
\label{def_gp}
Denote by $V(\Omega,P,f)$ the set of $\Omega$-tropical series $g$
such that $g\geq f$ and $g$ is not smooth at each of the points $\p\in P.$
 For a finite subset $P$ of $\Omega^\circ$ we define an operator $G_P
$ acting on $\Omega$-tropical series as $$G_P f(z)=\inf \{g(z)|g\in V(\Omega,P,f)\}.$$ If $P$ contains only one point $\p$, we write $G_\p$ instead of $G_{\{\p\}}$. 
\end{definition}

\begin{figure}[h!]
\begin{tikzpicture}
[x={(0.5cm,0cm)}, y= {(0.1cm,0.25cm)}, z={(0.09cm,0.45cm)}, scale=6]
\draw (0,1,0)--(0,0,0)--(1,0,0)--(1,1,0);
\draw[dashed](1,1,0)--(0,1,0);
\draw[very thick] (0,0,0)--(1/3,1/3,1/3)--(2/3,1/3,1/3)--(1,0,0);
\draw[very thick] (0,1,0)--(1/3,2/3,1/3)--(2/3,2/3,1/3)--(1,1,0);
\draw[very thick] (1/3,1/3,1/3)--(1/3,2/3,1/3);
\draw[very thick] (2/3,1/3,1/3)--(2/3,2/3,1/3);
\end{tikzpicture}
\begin{tikzpicture}[scale=2]
\draw[very thick] (0,0)--(1/3,1/3)--(2/3,1/3)--(1,0);
\draw[very thick] (0,1)--(1/3,2/3)--(2/3,2/3)--(1,1);
\draw[very thick] (1/3,1/3)--(1/3,2/3);
\draw[very thick] (2/3,1/3)--(2/3,2/3);
\end{tikzpicture}
\begin{tikzpicture}
[x={(0.5cm,0cm)}, y= {(0.1cm,0.25cm)}, z={(0.09cm,0.45cm)}, scale=6]
\draw (0,1,0)--(0,0,0)--(1,0,0)--(1,1,0);
\draw[dashed](1,1,0)--(0,1,0);
\draw[very thick] (0,0,0)--(2/15,4/15,4/15)--(1/5,1/3,1/3)--(2/3,1/3,1/3)--(1,0,0);
\draw[very thick] (0,1,0)--(2/15,11/15,4/15)--(1/5,2/3,1/3)--(2/3,2/3,1/3)--(1,1,0);
\draw[very thick] (1/5,1/3,1/3)--(1/5,2/3,1/3);
\draw[very thick] (2/3,1/3,1/3)--(2/3,2/3,1/3);
\draw[very thick] (2/15,4/15,4/15)--(2/15,11/15,4/15);
\end{tikzpicture}
\begin{tikzpicture}[scale=2]
\draw[very thick] (1,0)--(2/3,1/3)--(1/5,1/3)--(2/15,4/15)--(0,0);
\draw[very thick] (1,1)--(2/3,2/3)--(1/5,2/3)--(2/15,11/15)--(0,1);
\draw[very thick] (2/3,1/3)--(2/3,2/3);
\draw[very thick] (1/5,1/3)--(1/5,2/3);
\draw[very thick] (2/15,4/15)--(2/15,11/15);
\draw (1/5,1/2) node {$\bullet$};
\end{tikzpicture}
\caption{On the left: $\Omega$-tropical series $f=\min(x,y,1-x,1-y,1/3)$ and the corresponding tropical curve. On the right:  the fat point is $p=(\frac{1}{5},\frac{1}{2})$, the result $G_{p}f=\min(x + \frac{2}{15},y,1-x,1-y,\frac{1}{3})$ of the application of $G_p$ to $f$ and the corresponding tropical curve, note that this curve passes through $p$.}
\label{fig_3dpicture}
\end{figure}


\begin{proposition}[\cite{us_series}]\label{prop_gpsconvergence}
Let $\q_1,\q_2,\dots$ be an infinite sequence of points in $P$ where each point $\p_i,i=1,\dots, n$ appears infinite number of times. Let $f$ be any $\Omega$-tropical series. Consider a sequence of $\Omega$-tropical series $\{ f_m\}_{m=1}^\infty$ defined recursively as $$f_1=f, f_{m+1}=G_{{\bf q}_m} f_m.$$
Then the sequence $\{ f_m \}_{m=1}^\infty$ uniformly (on $\Omega^\circ$) converges to $G_P f$.
\end{proposition}

\begin{corollary}
Note that if $G_{{\bf q}_n}\dots G_{{\bf q}_1}f$ is close to the limit $G_Pf$, then Lemma~\ref{lemma_ecloseseries} implies that the corresponding tropical curves are also close to each other. Recall that we always present tropical series in the canonical form (i.e. all coefficients are taken to be as small as possible) and the closeness of the series as functions implies the closeness of their corresponding coefficients. 
\end{corollary}


\begin{definition}
\label{def_add}
Consider an $\Omega$-tropical series $f$ in the canonical form \eqref{eq_series}, fix $(k,l)\in\A$ and $c\geq 0$. We denote by $\Add_{kl}^c f$ the $\Omega$-tropical series $$(\Add_{kl}^c f) (x,y)=\min\left(a_{kl}+c+kx+ly,\inf\limits_{\substack{{(i,j)\in\A} \\ {(i,j)\ne (k,l)}}}(a_{ij}+ix+jy)\right).$$		  
\end{definition}

\begin{lemma}[\cite{us_series}]
\label{lem_singlegp}
Let $f = \inf_{(i,j)\in\A_\Omega}(ix+jy+a_{ij})$ be an $\Omega$-tropical series in the canonical form, such that the curve $C(f)$ doesn't pass through a point $\p=(x_0,y_0)\in\Omega^\circ$, and $f$ is $a_{kl}+kx+ly$ near $\p$. Then 
$G_\p f$ differs from $f$ only in a single coefficient $a_{kl}$, i.e. $G_\p f = \Add_{kl}^c f$ for some $c>0$.
\end{lemma}

\begin{remark}
Note that if $\Omega$ is a lattice polygon, the points $P$ are lattice points, all the coefficients $a_{ij}$ in $f$ are integer numbers, then, in the setup of Proposition~\ref{prop_gpsconvergence},  all the increments $c$ of the coefficients in $G_{q_{m+1}} f_m=\Add_{kl}^cf_m$ are integers, and therefore the sequence $\{f_m\}$ always stabilizes after a {\bf finite} number of steps.  
\end{remark}

\begin{figure}[h]
    \centering
\begin{tikzpicture}[scale=0.4]
\draw[very thick](0,0)--++(1,1)--++(0,7)--++(4,0)--++(1,-1)--++(0,-6)--++(-5,0);
\draw[very thick](1,8)--++(-1,1);
\draw[very thick](5,8)--++(0,1);
\draw[very thick](6,7)--++(1,0);
\draw[very thick](6,1)--++(1,-1);
\draw(3,4)node{$\bullet$};
\draw(3,4)node[right]{$p$};
\draw(4.5,6.5)node{\Large$\Phi$};
\draw[->][very thick](8,4)--(10,4);

\begin{scope}[xshift=300]
\draw[very thick](0,0)--++(2,2)--++(0,5)--++(3,0)--++(0,-5)--++(-3,0);
\draw[very thick](2,7)--++(-2,2);
\draw[very thick](5,7)--++(0,2);
\draw[very thick](5,7)--++(2,0);
\draw[very thick](5,2)--++(2,-2);
\draw(3,4)node{$\bullet$};
\draw(3,4)node[right]{$p$};
\draw[->][very thick](8,4)--(10,4);
\end{scope}

\begin{scope}[xshift=600]
\draw[very thick](0,0)--++(3,3)--++(0,3)--++(1,0)--++(0,-3)--++(-1,0);
\draw[very thick](3,6)--++(-3,3);
\draw[very thick](5,7)--++(0,2);
\draw[very thick](5,7)--++(2,0);
\draw[very thick](5,7)--++(-1,-1);
\draw[very thick](4,3)--++(3,-3);
\draw(3,4)node{$\bullet$};
\draw(3,4)node[left]{$p$};
\end{scope}

\end{tikzpicture}
\caption{Illustration for Remark~\ref{rem_smooth}. The operator $G_p$ shrinks the face $\Phi$ where $p$ belongs to. First, $t=0$, then $t=0.5$, and finally $t=1$ in $\Add_{ij}^{ct}f$. Note that combinatorics of the new curve can change when $t$ goes from $0$ to $1$.} 
\label{fig_ShrinkPhi}
\end{figure}

\begin{remark}
\label{rem_smooth} Let $i,j,c$ be such that $G_\p f=\Add_{ij}^c f$.
It is usefull to think of the operator $\Add_{ij}^c$  as a member of a continuous family of operators $$f\to \Add_{ij}^{ct}f,\text{\ where $t\in[0,1]$}.$$ This allows us to observe the tropical curve {\it during} the application of $\Add_{ij}^c$, in other words, we look at the family of curves defined by tropical series $\Add_{ij}^{ct}f$ for $t\in[0,1]$. See Figure~\ref{fig_ShrinkPhi}. Note that a curve on the middle picture has a nodal point (at the four-valent vertex) -- this is the worst thing that can happen with a non-terminal ($t\neq 1$) member of the family if the initial curve was smooth \cite{us_series}. This fact (Proposition \ref{prop_atmostnodal}) is crucial for deriving the lower estimate on the toppling function. \end{remark}

\subsection{$\Pi$-tropical polynomials}
Let $\Pi$ be a $\mathbb{Q}$-polygon. Similarly to the case of tropical series, a $\Pi$-tropical polynomial $g$ is a tropical polynomial vanishing along $\partial\Pi.$ A natural invariant of $g$ measuring its complexity is its quasi-degree.

\begin{definition}[Definition 9.3 in \cite{us_series}]\label{def_quasidegree}
Denote by $S(\Pi)$ the set of all sides of $\Pi.$ The quasi-degree of $\Pi$-torpical polynomial $g$ is a function $\operatorname{qd}_g\colon S(\Pi)\rightarrow\mathbb{Z}_{\geq 0}$ such that $g$ near a side $s\in S(\Pi)$ coincides with $\operatorname{qd}_g(s)\lambda_s,$ where $\lambda_s$ is a support function of $\Pi$ vanishing along the side $s$ and having primitive gradient. 
\end{definition}

Observe that the symplectic area (where each Euclidean length of an edge is taken with a weight equal to the primitive vector parallel to it) of the tropical curve given by $g$ depends only on the quasi-degree of $g,$ being a weighted sum of lengths of the sides of $\Pi.$ Therefore, its minimization for $C(f_{\Pi,P})$ in the class of curves passing through $P$ follows from the pointwise minimality of $f_{\Pi,P}$.   

Being nice is a technical property for a quasi-degree needed for constructing an auxiliary state in the construction of the lower estimate.
\begin{definition}[Definition 9.4 in \cite{us_series}]\label{def_nice}
A function $d\colon S(\Pi)\rightarrow \mathbb{Z}_{>0}$ is called {\it nice} if for each side $s\in S(\Pi)$ with $d(s)>1$ the values of $d$ on the two sides of $\Pi$ adjacent to $s$ are equal to $1.$ 
\end{definition}

The following guarantees that on unimodular polygons there exist tropical polynomials, with any prescribed nice quasi-degree, defining a smooth tropical curve situated in an arbitrarily small neighborhood of the boundary. We use this curve to construct the auxiliary state by which we replace the maximal stable state.
\begin{proposition}[Theorem 9.5 in \cite{us_series}]\label{prop_smoothnice}
For a unimodular $\Pi,$ nice $d\colon S(\Pi)\rightarrow \mathbb{Z}_{>0}$ and $\varepsilon>0$ there exist a $\Pi$-tropical polynomial $g_\varepsilon<\varepsilon$ such that $\operatorname{qd}_{g_\varepsilon}=d,$ the curve $C(g)$ is smooth and is contained in $B_\varepsilon(\partial\Pi).$
\end{proposition}

By performing small corner cuts (blow-ups in toric terminology) of $\Pi$ we can replace it with an unimodular polygon on which a given partial tropical relaxation has a nice quasidegree.  
\begin{proposition}[Proposition 4 in \cite{us_series}] \label{prop_niceapprox}
Let $\Pi$ be a $\mathbb{Q}$-polygon. Consider a sequence of not necesarily distinct points $q_1,q_2,\dots, q_m\in\Pi^\circ.$ Then, for all small enough $\varepsilon>0,$ there exist a unimodular $\mathbb{Q}$-polygon $\Pi'$ such that all $q_j$ are in $(\Pi')^\circ$ and for the operators $G=G_{q_m}\dots G_{q_1}$ considered on both $\Pi$-tropical and $\Pi'$-tropical polynomials the following holds:
\begin{itemize}
\item  $\Pi'$-tropical polynomial $G 0_{\Pi'}$ has a nice quasi-degree;
\item $0\leq G 0_\Pi-G 0_{\Pi'}<\varepsilon;$
\item $G 0_\Pi\leq \varepsilon$ on $\Pi\backslash \Pi'.$
\end{itemize} 
\end{proposition}

The following proposition ensures that for the lower estimate on the toppling function, we may stay within the space of at most nodal tropical curves (for which their sandpile counterparts made of solitons are well behaved under the action of waves).

\begin{proposition}[Lemma 8.6 in \cite{us_series}]\label{prop_atmostnodal}
Let $\Phi$ be a face of a tropical curve $C(g)$ containing a point $p.$ Let $G_p g=\Add_{i,j}^c g.$ For $t\in[0,1),$ denote by $\Phi(t)$ the face of $C(\Add_{i,j}^{ct} g)$ containing $p.$ If $\Phi$ is unimodular, then for all $t\in (0,1)$ all vertices of $\Phi(t)$ are either smooth or nodal vertices of the tropical curve $C(\Add_{i,j}^{ct} g).$
\end{proposition} 
Here, by a nodal vertex we mean a simple double point, i.e. such that the defining polynomial around this vertex is equal to $\min{0,x,y,x+z}$ up to a constant and a modular transformation. A smooth vertex of a tropical curve is the one locally given by $\min(0,x,y).$ Note that every face of a tropical curve with only smooth or nodal vertices is unimodular.

\section{Smoothing of superharmonic functions}\label{sec_smoothings}

Let $F:\ZZ^2\to\ZZ$ be a superharmonic (i.e. $\Delta F\leq 0$ everywhere) function.  Call $D(F)=\{v\in \ZZ^2|\Delta F (v)<0\}$ the deviation set of $F$.

\begin{definition}[\cite{us_solitons}]
\label{def_thetan}
For each $n\in\mathbb{N}$ consider the set $\Theta_n(F)$ of all integer-valued superharmonic functions $G$ such that  $F-n\leq G\leq F$ and $G$ coincides with $F$ outside a finite neighborhood of the deviation set of $D$, i.e. 
$$\Theta_n(F)=\{G:\ZZ^2\to \ZZ| \Delta G \leq 0, F-n\leq G\leq F, \exists c>0, \{F\ne G\}\subset B_{c}(D(F))\}.$$ Define $S_n(F):\ZZ^2\to\ZZ$ to be $$S_n(F)(v)=\min\{G(v)|{G\in\Theta_n(F)}\}.$$ We call $S_n(F)$ {\it the $n$-smoothing of $F$}. 
\end{definition}

In \cite{Kalinin_2023} the same procedure is called {\it husking} because this name suits better to this operation (note that $S_n(F)$ is not smooth, it is a function defined on a lattice). However here we keep old terminology, as in \cite{us_solitons}.

Let us fix $p_1,p_2,q_1,q_2,c_1,c_2 \in \ZZ$ such that $p_1q_2-p_2q_1=1$. Consider the following functions on $\ZZ^2$:  
\begin{equation}\label{eq_edge}
\Psi_{\edge}(x,y)=\min(0,p_1x+q_1y),
\end{equation}
\begin{equation}\label{eq_tripode}
\Psi_{\vertex}(x,y)=\min(0,p_1x+q_1y,p_2x+q_2y+c_1),
\end{equation}
\begin{equation}\label{eq_node}
\Psi_{\node}(x,y)=\min\Big(0, p_1x+q_1y,p_2x+q_2y+c_1, (p_1+p_2)x+(q_1+q_2)y+c_2 \Big).
\end{equation}

\begin{theorem}[\cite{us_solitons}]\label{th_stabilfn}
Let $F$ be a) $\Psi_{\edge}$, b)$\Psi_{\vertex}$, or c)$\Psi_{\node}$. The sequence of $n$-smoothings $S_n(F)$ of $F$ stabilizes eventually as $n\to\infty$, i.e. there exists $N>0$ such that $S_n(F)\equiv S_N(F)$ for all $n>N$. In other words, 
there exists a pointwise minimum $\theta_F=S_N(F)$ (we call it {\it the canonical smoothing of $F$}) in $\bigcup\Theta_n(F)$. 
\end{theorem}

\begin{figure}
\label{fig_triads}
\includegraphics[width=0.3\textwidth]{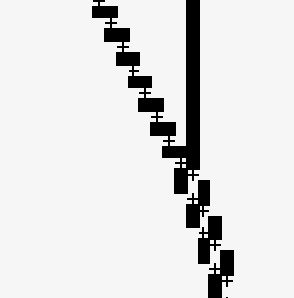}
\includegraphics[width=0.49\textwidth]{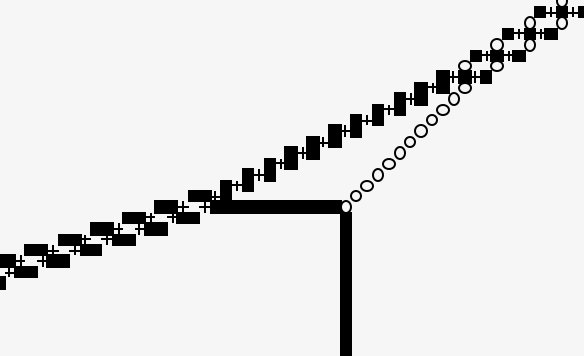}
\caption{Examples of triads. One triand is on the left and three triads are on the right.}
\end{figure}

\subsection{Waves and operators $G_\p$.}
\begin{remark}Note that $\Delta \theta_F\geq -3$ because otherwise we could decrease $\theta_F$ at a point violating this condition, preserving superharmonicity of $\theta_F$, and this would contradict to the minimality of $\theta_F$ in $\bigcup\Theta_n(F)$.
\end{remark}

Consider a state $\phi=\state{3} + \Delta \theta_F$. By the previous remark, $\phi\geq 0$. Let $v\in \ZZ^2$ be a point far from $\{\Delta\theta_F\ne 0\}$. Let $F$ equal to $ix+jy+a_{ij}$ near $v$. Then, informally, sending a wave from $v$ increases the coefficient $a_{ij}$ by one. 
\begin{lemma}[\cite{us_solitons}]
\label{lemma_wavegp}
In the above conditions, $W_v\phi=\state{3}+\theta_{F'}$ where $W_v$ is sending a wave from $v$ (see \cite{us_solitons}) and $F'=\Add_{ij}^1F$ (Definition~\ref{def_add}).
\end{lemma}

Canonical smoothings in Theorem~\ref{th_stabilfn} provide us with the building blocks of the deviation locus $D(\phi_\h^\circ)$ in case when $C(f_{\Omega,P})$ is at most nodal. The smoothing of $\Psi_{\edge}$ represents {\it a sandpile soliton}, becoming a tropical edge in the limit, the smoothing of $\Psi_{\vertex}$ represents {\it a sandpile triad}, becoming a smooth tropical vertex in the limit, and the smoothing of $\Psi_{\node}$ represents the degeneration of two sandpile triads into the union of two sandpile solitons. 


\begin{definition}
A vertex $V$ of a tropical curve $C(f)$ is {\it smooth} if the restriction of $f$ to a small neighborhood of $V$ can be presented as $\min(x,y,0)$ after an $SL(2,\ZZ)$ change of coordinates following by a translation. A vertex $V$ of $C(f)$ is called a {\it node} if the restriction of $f$ to a small neighborhood of $V$ can be presented as $\min(x, y, 0, x + y)$ after an $SL(2,\ZZ)$ change of coordinates following by a translation. A tropical curve is called {\it at most nodal} if all its vertices are smooth or nodal.
\end{definition}

\begin{definition}[Canonical smoothing of a tropical polynomial]
\label{def_functionsmoothing}
Let $f$ be a $\Pi$-tropical polynomial, such that $C(f)$ is at most nodal tropical curve. Then, for small enough $\h>0$,  we define {\it the canonical smoothing} $\smooth_\h(f):\Pi_\h\to\ZZ$ as follows. We consider the lattice $\h\ZZ^2\subset\RR^2$ and define $F(x,y) = [\h^{-1}f(x,y)], (x,y)\in\Pi$. This defines a piece-wise linear function on $\Pi$ (cf. \eqref{eq_rounded}). The curve $C(f)$ has a finite number of edges and vertices, which are smooth or nodal. Hence the same is true for $C(F)$ if $\h$ is small enough. Each local equation of vertices and edges of $C(F)$ belongs to the cases in Theorem~\ref{th_stabilfn}. We apply Theorem~\ref{th_stabilfn} for all these local equations. Hence there exists $N>0$ such that the smoothings of all the local equations of edges and vertices of $C(F)$ stabilize after $N$ steps. Finally, we define $\smooth_\h(f)$ as $S_N(F)$ restricted to $\Pi_\h$ (here we treat $F$ as a function on $\h \ZZ^2\to\ZZ$). We call this procedure the {\it canonical smoothing} of $F$ (or the canonical smoothing of $f$ with respect to $\h$). Note that $\smooth_\h(f)$ may have negative values.
\end{definition}

\begin{remark}
Note that $\Delta\smooth_\h(f)$ is zero almost everywhere and its non-zero locus consists of patches of types $\Delta\Psi_{\edge}, \Delta\Psi_{\vertex}, \Delta\Psi_{\node}$ (in the notation of Theorem~\ref{th_stabilfn}). Also, we should not worry about how these patches glue together because of the invariant definition of $\smooth_\h(f)$, it is the function obtained after $N$ steps of the smoothing process (Definition~\ref{def_thetan}). Another way to understand it is to look at the proof of Theorem~\ref{th_stabilfn} (\cite{us_solitons}, Section 9): in order to find the smoothing of a vertex we first smooth the function along edges, and then prove that smoothing of the vertex affects only a finite neighborhood of this vertex.
\end{remark}

\begin{remark}
\label{rem_constantsmoothing}
The canonical smoothing procedure doesn't change $F$ outside $B_{N\h}(D(F))$ where $N$ is an absolute constant, depending on the slopes of the edges of $C(f)$. Therefore if $\h$ is small enough, the smoothings of different vertices and edges never overlap. In this case, we say that $\smooth_\h(f)$ is well defined.
\end{remark}

%

\section{Proof of the lower estimate}
\label{sec_lower}


Fix $\e>0$. Let $\Pi$ be a $\QQ$-polygon and $P\subset \Pi^\circ$ a finite collection of points. Using Proposition~\ref{prop_gpsconvergence} choose $f_m$ which is $\e/2$-close to $f_{\Pi,P}$ and their quasidegree (see Definition \ref{def_quasidegree}) coincide. Then, $f_m=G_m 0_\Pi$ where $G:=G_{{\bf q}_m}G_{{\bf q}_{m-1}}\dots G_{{\bf q}_1}$ for some ${\bf q}_i\in P$, $0_\Pi$ is the function which is zero on $\Pi$.  Note that by Lemma \ref{lem_singlegp} 
\begin{equation}
\label{eq_Gp}G_{{\bf q}_m}G_{{\bf q}_{m-1}}\dots G_{{\bf q}_1} 0_\Pi=\Add_{i_mj_m}^{e_m}\dots \Add_{i_1j_1}^{e_1} 0_\Pi
\end{equation}
 for certain $(i_k,j_k)\in \ZZ^2, e_i\geq 0$.  
\begin{lemma}
\label{lemma_fm}
 There exist a $\QQ$-polygon $\Pi'\subset \Pi$ and  a $\Pi'$-tropical polynomial $g$ such that $g<\e$ and $C(g)$ is a smooth tropical curve in $\Pi'$; and for all $\eta>0$ small enough the following function
$$\tilde G g, \text{\ where\ }\tilde G:=\Add_{i_mj_m}^{e_m-\eta}\dots \Add_{i_1j_1}^{e_1-\eta}$$
is $\e$-close to $f_{\Pi,P}$ on $\Pi$ (in particular, $\tilde Gg<\e$ on  $\Pi\setminus\Pi'$), and during the application of all operators $\Add$ (Remark~\ref{rem_smooth}) all the corresponding tropical curves are at most nodal (see Proposition \ref{prop_atmostnodal}).
\end{lemma}

\begin{proof}

Using Proposition \ref{prop_niceapprox} for $G$ we construct a polygon $\Pi'\subset\Pi$ with the property that $G0_{\Pi'}$ is nice (Definition \ref{def_nice}). Then, under this assumption, Proposition \ref{prop_smoothnice} asserts existence of a $\Pi'$-tropical polynomial $g$ of the same quasidegree as $G0_{\Pi'}$, $g|_{\Pi'}<\e$, the curve $C(g)$ is smooth, $Gg=G0_{\Pi'}$ on $\Pi'$, and
$G g$ is $\e/2$-close to $G 0_{\Pi}$ on $\Pi$.  Since the quasi-degree of $g$ and $Gg$ coincide, during the calculation of $Gg$ we never apply a wave operator for a face that has a common side with $\partial\Pi'$, so $Gg=g<\e$ on $\Pi\setminus\Pi'$. 
Note that in the product $G_{{\bf q}_m}G_{{\bf q}_{m-1}}\dots G_{{\bf q}_1}g$ each $G_{{\bf q}_k}$ is the application of $\Add_{i_k,j_k}^{e_{k}}$ for some $e_k>0$ (we assume that the action is effective, if $e_k=0$ we remove the corresponding $\Add_{i_k,j_k}^{e_{k}}$ from the product), i.e. we increase the coefficient in the monomial $i_kx+j_ky$ by $e_k$.

%
Fix $\eta>0$ small enough. We
replace in  \eqref{eq_Gp} $$\Add_{i_kj_k}^{e_k} \text{\ by\ }\Add_{i_kj_k}^{e_k-\eta},  \text{\ for\ } k=1,\dots, m.$$ Denote $\tilde f_0=g, \tilde f_{k+1}=\Add_{i_kj_k}^{e_k-\eta}\tilde f_k$. 
  Then, all the tropical curves defined by $\tilde f_k,k=1,\dots, m$ are at most nodal on $\Pi$ as well as each tropical curve in the family during the application of $\Add_{i_kj_k}^{e_k-\eta}$ to $\tilde f_k$ by Proposition \ref{prop_atmostnodal}.  We briefly recall the idea: the only two possibilities for how the tropical curve can fail to be at most nodal during our procedure in \eqref{eq_Gp} is the appearance of a non-smooth vertex inside $\Pi^\circ$ and the appearance of an edge with multiplicity bigger than one inside $\Pi^\circ$ or at the corners of $\Pi$; both things can happen at the very last moment of application of $\Add_{i_kj_k}^{e_k}$, thus diminishing $e_k$ by any positive number prevents appearing non-smooth vertices and edges of multiplicity bigger than one. 
  
Finally, for $\eta$ small enough the tropical curve defined by $\tilde f_m$ is $\e/2$-close to the tropical curve defined by $f_m=G_{{\bf q}_m}G_{{\bf q}_{m-1}}\dots G_{{\bf q}_1}g$ and thus $\e$-close to  $f_{\Pi,P}$.


\end{proof}

\begin{proof}[Proof of Proposition~\ref{prop_qpolygonconvergence}]
Choose any $\e>0$ small enough. We will construct an auxiliary state $\psi_\h$ on $\Pi_\h$ whose toppling function is less than that of $\phi_\h$ by at most $\e$, and whose partial relaxation (via wave decomposition
) 
is completely controlled by operators $\Add_{i_k,j_k}^{e_{k}-\eta}$. 

We use the notation of Lemma~\ref{lemma_fm} which can be applied in our situation. Let $\psi_\h$ be $\state{3}+\sum_{p\in P^{\h}} \delta_{p}+\Delta\smooth_\h(g)$ (Definition~\ref{def_functionsmoothing}) for $g$ from Lemma~\ref{lemma_fm}. Note that $\psi_\h\leq \phi_\h$. Furthermore, for $\h$ small enough, $\psi_\h$ consists of three grains almost everywhere and is made of sandpile solitons and triads near $C(g)$. Relaxing $\phi_\h$ consists of applying wave operators from the points in $P^\h$. Thus, performing an even smaller number of waves from $P^\h$ to $\psi_\h$ and counting topplings furnishes us a lower bound for the toppling function of $\phi_\h$. 

Namely, using $\tilde f_k$ we define $F_k =[\h^{-1}\tilde f_k]:\Pi_\h\to\ZZ_{\geq 0}$ as in \eqref{eq_rounded}.
We can choose $\h>0$ small enough such that all canonical smoothings $\smooth_\h(f_k),k=0,\dots,m$  are well defined (Remark~\ref{rem_constantsmoothing}). We denote by $\smooth_\h(\tilde f_k)$ the canonical smoothing of $F_k$ (see Definition~\ref{def_functionsmoothing}). Define the states $\psi_k=\state{3}+\Delta \smooth_\h(\tilde f_k).$ 

The final step is to use the fact that wave operators commute with smoothings (we proceed as in Lemma~\ref{lemma_wavegp} but with the global notation). Namely, let $0\leq k\leq m$. Fix the notation by $$F_k(x,y)=\min_{(i,j)\in\A}(ix\h^{-1}+jy\h^{-1} + [a_{ij}\h^{-1}])$$ with finite $\A\subset\mathbb{Z}^2$. Let $\eta=M\h$ in Lemma~\ref{lemma_fm} for $M\in \ZZ_{>0}$ big enough (recall that $\h\to 0$.) All the points $p_i$ do not belong to $C(\tilde f_k)$ (because $M$ is big enough, and we did $M$ waves less than required to get $p_i$ on the tropical curve, i.e. it is enough to have $M$ bigger than the width of all smoothings of tropical edges) and so do not belong to $C(F_k)$ as long as $\h$ is small enough. Let $v_k=\h[\h^{-1}{\bf q}_k]$. Then $\psi_k(v_k)=3$. Suppose that $v_k$ belongs to the region where $i_0x\h^{-1}+j_0y\h^{-1} + [a_{i_0j_0}\h^{-1}]$ is the minimal tropical monomial. Denote $$F'_k(x,y)=\min_{(i,j)\in\A}(ix\h^{-1}+jy\h^{-1} + [a_{ij}'\h^{-1}])$$ where $a_{ij}'=a_{ij}$ if $(i,j)\ne(i_0,j_0)$ and $a_{i_0j_0}'=a_{i_0j_0}+\h$. Then, by Lemma \ref{lemma_wavegp} 
 $$W_{v_k}\psi_k = \state{3}+\Delta\smooth_\h(F'_k)$$. 

Therefore $$\state{3}+\Delta \smooth_\h(F_{k+1})= W_{v_{k+1}}^{[e_{k+1}\h^{-1}]-M} \big(\state{3}+ \Delta\smooth_h(F_k)\big).$$ 

Hence the toppling function of $\phi_h=\state{3}+\sum_{p\in P} \delta_{p^\h}$ is at least $[\h^{-1}(\tilde f_m)-g]$. The proposition follows since the toppling function of $\psi_h$ is less than the toppling function of $\phi_h$ and $\tilde f_m$ is close to $f_{\Pi,P}.$ Note that thanks to the construction of $g$ we had no topplings near the boundary of $\Pi'$ during this process. We finished the proof of Proposition~\ref{prop_qpolygonconvergence}.
\end{proof}


Let us informally summarise the proof. If $\phi=\state{3}+\Delta F$ where $F$ is a tropical polynomial such that $C(F)$ is smooth or nodal, then sending a wave from $\p$ corresponds to the operator $\Add_{i,j}^\h$ whre $(i,j)$ is the gradient of $F$ at $\p$. Then, $(\phi+\delta_\p)^\circ$ can be obtained by sending waves until the value at $\p$ is less than three. This corresponds to applying $\Add_{i,j}\h$ until the tropical curve reaches the point $\p$, i.e. this corresponds to the operator $G_\p$. Thus lower bound for the toppling function may be obtained by starting from an auxiliary state $\psi$, composed of sandpile solitons and triads, which behave under the action of wave in a known way, thanks to \cite{us_solitons}, and then following the dynamic prescribed by tropical operators $G_\p$. In this process, we lower the required (rescaled) number of waves from each point from $e_k$ to $e_k-M\h$ in order to keep our tropical curves at most nodal, because only for such curves we do know their sandpile approximations that behave in a controlled way under the wave action.


\section{The multiplicities of the edges via a weak convergence.}
\label{sec_weakconvergence}
In the notation of Theorem~\ref{th_main}, define $$\tilde\phi_\h(x,y)=\h^{-1}\left(3-\phi_\h^\circ(\h[\h^{-1}x],\h[\h^{-1}y])\right).$$ Note that $\tilde\phi_\h$ is not zero only along $D(\phi_\h^\circ)$. 

\begin{theorem}[A stronger form of Theorem 2 announced in \cite{announce}] There exists a $*$-weak limit $\tilde\phi$ of the sequence $\tilde\phi_\h$ as $\h\to 0$. Moreover, there exists a unique assignment of multiplicities $m_e$ for the
edges $e$ of $C(f_{\Omega,P})$ such that for all smooth functions $\Phi$ supported on $\Omega^\circ$ we have
$$\tilde\phi(\Phi)=\lim\limits_{\h\to 0}\int\limits_{\RR^2}\tilde\phi_\h\Phi=\sum\limits_{e\in E}\left(||l_e||\cdot m_e\cdot\int_e\Phi\right),$$
where $E$ is the set of all edges of $C(f_{\Omega,P})$ and $l_e$ is a primitive vector of $e\in E$, i.e. the coordinates of $l_e$ are coprime integers and $l_e$ is parallel
to $e$.
\end{theorem}
\begin{proof}Choose small $\e>0$, the same as in the proof of Theorem~\ref{th_main}. Outside of $C(f_{\Omega,P})$ the $*$-weak limit of $\tilde\phi_\h$ is zero because for $\h$ small enough $\tilde\phi_\h\equiv 3$ outside of $\e$-neighborhood of $C(f_{\Omega,P})$ by Theorem~\ref{th_main}. 

Suppose for now that $C(f_{\Omega,P})$ is a smooth tropical curve (in particular, for each edge $e$ its multiplicity $m_e$ is one). Consider an edge $e$ of $C(f_{\Omega, P})$ and a strict subinterval $e'$ of it, $e'\subset e$. Consider a thin rectangular $Q$ whose one side is parallel to $e$ and another side has length $C\h$ with a constant $C_1$ big enough such that $Q\cap C(f_{\Omega,P})=e'$. Choosing $\h$ small enough and using Lemma \ref{lem_integraloflaplacian} below, we see that $\big|\int_Q\psi_\h-\h^{-1}|l_e|\cdot |e'|\big|<C_2\e$, because only the contribution over the long sides of $Q$ matters and the contribution of infinitesimally small sides of $Q$ is small, and the actual toppling function rescaled by $\h$ is $\e$-close to $f_{\Omega,P}$. 

If $C(f_{\Omega,P})$ is not a smooth tropical curve, then we approximate $f_{\Omega,P}$ by tropical series $f_{\Omega,P}^\e, e\to 0$ which give smooth tropical curves, as in the proof of Lemma \ref{lemma_fm}. Then each edge $e$ of $C(f_{\Omega,P})$ of multiplicity $m_e$ appears as the limit of $m_e$ edges of $C(f_{\Omega,P}^\e)$ of the same direction. As above we consider a thin rectangular, now containing all these $m_e$ edges of $C(f_{\Omega,P}^\e)$, and the same argument concludes the proof.
\end{proof}

\begin{lemma}[Lemma~6.3 of \cite{us_solitons}]
\label{lem_integraloflaplacian}
Let $A$ be a finite subset of $\mathbb{Z}^2$ and $\partial A$ be its subset consisting of points adjacent to $\mathbb{Z}^2\backslash A.$ Then, for any function $F\colon A\rightarrow\mathbb{Z}$
$$\sum_{v\in A\backslash\partial A}\Delta F(v)=\sum_{\substack{v\in\partial A\\v'\in\partial A, v\sim v'}}(F(v)-F(v')).$$
\end{lemma}

\begin{proof} The proof follows from the definition of the discrete Laplacian -- all contributions of $F(v)$ for $v\in A\partial A'$ cancel out and we are left with a summation over the boundary.
\end{proof}

%

\section{Discussion}
\label{sec_discussion}


\subsection{Continuous limit} Tropical curves appear as limits of algebraic curves under the map $\log_t |\cdot|\colon({\mathbb C}^*)^2\rightarrow\mathbb{R}^2$ when $t\to\infty$. It is natural to ask how we can obtain a continuous family of ``sandpile'' models that converges to the pictures we studied in this paper. We do not know how to do this. Meanwhile, recently amoebas showed up in the study of the leaky sandpile model \cite{alevy2020limit}. Again, we do not know how to relate their results to ours. Going to the limit in the other direction, by considering sandpile dynamics on the set of tropical curves we (experimentally) found a power law for the avalanche statistics in the continuous tropical sandpile model \cite{sandcomputation}, no proofs\footnote{There is some proof for a power law in  the one dimensional case \cite{shkolnikov2023relaxation}.} are available for now, only computational statistics about that model \cite{stattrop}.

\subsection{Fractals and patterns in sandpile}
The sandpile on $\ZZ^2$ exhibits a fractal structure; see, for example, the pictures of the identity element in the critical group \cite{LP}. As far as we know, only a few cases have a rigorous explanation. It was first observed in \cite{ostojic2003patterns} that if we rescale by $\sqrt n$ the result of the relaxation of the state with $n$ grains at $(0,0)$ and zero grains elsewhere, it weakly converges as $n\to\infty$. Then this was studied in \cite{levine2009strong} and was finally proven in \cite{PS}. However, the fractal-like pieces of the limit found their explanation later, in \cite{levine2013apollonian, LPS}, and happen to be curiously related to Apollonian circle packing. This allows to describe the identity on certain ellipses quite explicitly \cite{melchionna2020sandpile}.

In \cite{lang2019harmonic,lang2019sandpile} it was shown that tropical curves constitute ``linear'' flows on the sandpile group, while quadratic patches may be found in the ``quadratic'' flows.

A lot of work is to be done in the future. The work \cite{sadhu2017emergence}(see also \cite{sadhu2012pattern}) contains a lot of pictures and examples with an apparent piece-wise linear behavior.  We expect that the methods of this article will be used to study the fractal structure in those cases.

\subsection{The content of this paper and where to find proofs of previously announced results}

In \cite{announce} we announced several theorems in the case of lattice polygons and integral perturbation points, their generalized versions are proven in this paper. Here we list where to look for the proofs. Theorem 1, in \cite{announce},  is Theorem~\ref{th_main} here. Theorem 2 in \cite{announce} is proven in Section~\ref{sec_weakconvergence}.  Theorem 3 in \cite{announce} easily follows from Theorem~\ref{th_main}, and is proven in \cite{us_series}. Theorem 4 in \cite{announce} follows from Theorem~\ref{th_main} here just because of the definition of the function $f_{\Omega,P}$ (see Definition~\ref{def_fOmegaP}). We prove Theorem 5 in \cite{announce}  on the way of the proof of Theorem~\ref{th_main}, see Remark~\ref{rem_th5}. 

\subsection{Acknowledgments}

We thank Andrea Sportiello for sharing his insights on perturbative
regimes of the Abelian sandpile model which was the starting
point of our work. We also thank Grigory Mikhalkin, who encouraged us and gave a lot of advice about this paper. 

Also, we thank Misha Khristoforov and Sergey
Lanzat who participated in the initial state of this project, when we
had nothing except the computer simulation and pictures. Ilia Zharkov,
Ilia Itenberg, Kristin Shaw, Max Karev, Lionel Levine, Ernesto Lupercio, Pavol \v Severa, Yulieth Prieto, and Michael Polyak asked us a lot of questions; not all of these questions found their answers here, but
we are going to treat them in subsequent papers. 



\bibliography{../../bibliography}

\bibliographystyle{abbrv}
\end{document}